\documentclass[10pt]{amsart}
\usepackage{graphicx}
\usepackage{color}
\definecolor{orange}{rgb}{1,.7,0}
\definecolor{violet}{rgb}{.5,0,.5}
\definecolor{dg}{rgb}{0,0.67,0}
\usepackage{amssymb}
\usepackage{amssymb,verbatim,url}
\usepackage[all]{xy}

\newcommand{\cd}{\! \cdot \!}

\newcommand{\C}{{\mathbb C}}
\newcommand{\Z}{{\mathbb Z}}
\newcommand{\R}{{\mathbb R}}

\newcommand{\F}{{\mathbb F}}

\newcommand{\bbP}{{\mathbb P}}
\newcommand{\Q}{{\mathbb Q}}

\newcommand{\Gal}{\mbox{Gal}}

\newcommand{\Fr}{\mbox{Fr}}
\newcommand{\Aut}{\mbox{Aut}}

\newcommand{\Spec}{\mbox{Spec}}

\newcommand{\ced}{\mbox{cd}}
\newcommand{\cmmt}[1]{}
\numberwithin{equation}{section}
\numberwithin{table}{section}
\numberwithin{figure}{section}
\newtheorem{Theorem}{Theorem}[section]
\newtheorem{Proposition}{Proposition}[section]

\allowdisplaybreaks
\title[Division polynomials with Galois group $SU_3(3).2 \cong G_2(2)$]{Division polynomials with Galois group $SU_3(3).2 \! \cong \! G_2(2)$}
\author{David P. Roberts}
\address{Division of Science and Mathematics, University of
  Minnesota Morris, Morris, MN 56267, USA}
\email{roberts@morris.umn.edu}

\begin{document}

\begin{abstract}   We use a rigidity argument to prove the 
existence of two related degree twenty-eight 
covers of the projective plane with
Galois group $SU_3(3).2 \cong G_2(2)$.  
Constructing corresponding two-parameter polynomials
directly from the 
defining group-theoretic data seems 
beyond feasablity.  Instead we provide
two independent constructions of these
polynomials, one from $3$-division points
on covers of the projective line studied by Deligne and Mostow, and one
from $2$-division points of genus three curves studied by Shioda.
We explain how one of the covers also arises
as a $2$-division polynomial for a family of 
$G_2$ motives in the classification of 
Dettweiler and Reiter.   We conclude by specializing
our two covers to get interesting
three-point covers and number fields
which would be hard to construct
directly. 
\end{abstract}
\maketitle

\tableofcontents

\section{Introduction}
\label{intro}   

     Suppose $Y$ is a variety over $\Q$ with bad reduction at 
 a set $S$ of primes.   For any prime $\ell$, there are 
 associated number fields coming from the mod $\ell$ cohomology
 of the topological space $Y(\C)$.  On the one hand, these number fields are interesting because their
 Galois groups tend to be Lie-type groups and their
 bad reduction is constrained to be within $S \cup \{\ell\}$. 
 On the other hand, defining polynomials for these
 number fields are often beyond computational
 reach, even for quite simple $Y$ and very small $\ell$.   In this paper,
 we work out some remarkable examples in this 
 framework, with our computations of defining polynomials
 being {\em ad hoc} and just within the limits of 
 feasibility.  
 
 \subsection{Section-by-section overview}
     Section~\ref{background} provides background on the theoretical context, presenting
 it as a generalization of the familiar passage from an elliptic curve to 
 one of its division polynomials.  It then gives
 information on the Lie-type group which plays the
 central role for us, namely $SU_3(3).2 \cong G_2(2)$. 
 Finally, the section reviews an earlier construction of a one-parameter polynomial 
 for this Galois group due to Malle and Matzat \cite[p.\ 412]{MalMat}.
 
     Section~\ref{rigidcovers} explains how a rigidity argument 
   gives two canonical degree twenty-eight covers
   of surfaces defined over $\Q$, each with Galois group 
  $SU_3(3).2 \cong G_2(2)$.    In our notation,
  these covers are 
  \begin{align*}
  \pi_1 : X_1 & \rightarrow  U_{3,1,1}, &
  \pi_2 : X_2 & \rightarrow  U_{3,2}. 
  \end{align*}
  The bases are respectively 
  $U_{3,1,1} = M_{0,5}/S_3$ and 
  $U_{3,2} = M_{0,5}/(S_3 \times S_2)$, these being moduli spaces of five partially
  distinguishable points
  in the projective line.   We explain how the covers
  are related by a cubic correspondance 
  deduced from an exceptional 
  isomorphism $U_{2,1,1,1} \cong U_{2,1,2}$ 
  first studied by Deligne and Mostow \cite[\S10]{DM}.  
  Standard methods, as illustrated in \cite{CoversM05}, 
  might let one construct the
  covers $\pi_i$ directly if certain curves had genus zero. 
  However  these methods are obstructed by the
  fact that these curves have positive genus.
  
      Sections~\ref{DM}, \ref{shioda}, and \ref{G2} concern varietal sources for our
  covers.   Section~\ref{DM} starts with two different two-parameter
   families of covers of the projective line considered by Deligne and Mostow 
   \cite{DM1}.    
  Via the group $SU_3(3).2$, these families of curves yield $\pi_1$ and $\pi_2$
  from $3$-division points.  We use the second family to 
  compute a defining
  polynomial $F_2(a,b,x)$ for $\pi_2$,
  and then transfer this knowledge to also obtain a polynomial
  $F_1(p,q,x)$ for $\pi_1$.  Section~\ref{shioda} starts with 
  a large
   family of genus three curves studied by Shioda \cite{Shi2}.  
  This family already has an explicit $2$-division polynomial
  $S(r_1,r_3,r_4,r_5,r_6,r_7,r_9,z)$ with Galois group 
  $Sp_6(2)$.   We find appropriate loci in the parameter space where
  the Galois group drops to the subgroup $G_2(2)$, 
  and thereby independently get alternative polynomials
  $S_1(p,q,z)$ and $S_2(a,b,z)$ for the two covers.  
  Section~\ref{G2} explains how $F_1(p,q,z)$ also 
  arises as the $2$-division polynomial
  of a family of motives with 
  motivic Galois group $G_2$ studied by Dettweiler and
  Reiter \cite{DR}.   Sections~\ref{DM}, \ref{shioda}, and \ref{G2}
  each close with subsection explicitly relating $L$-polynomials
  modulo the relevant prime $\ell$ to our division polynomials.  
   
     Section~\ref{threepoint} shifts the focus away from varietal sources and onto 
 specializations of our explicit polynomials.   Specializing to suitable lines,
 we get fourteen new degree twenty-eight 
 three-point covers with Galois
 group $SU_3(3).2 \cong G_2(2)$.   These covers all have positive
 genus, and it would be difficult to construct them directly by 
 the standard techniques of three-point covers.
 
     Section~\ref{numbfield} specializes to points, finding 376 different degree
  $28$ number fields with Galois group  $SU_3(3).2 \cong G_2(2)$ and
  discriminant of the form $2^a 3^b$.  Again it would be 
  difficult to construct these fields by techniques within
  algebraic number theory itself.   We show that a thorough
  analysis of ramification in these fields is possible, despite the relatively large
  degree, by presenting such an analysis of the
  field with the smallest discriminant.

  \subsection{Computer platforms.}  The bulk of the calculations for this paper were carried 
  out in {\em Mathematica} \cite{Wo}.    However most calculations with
  number fields were done in {\em Pari} \cite{Pari} while most calculations
  with $L$-functions were done in {\em Magma} \cite{Mag}.   

  Many of the statements in this paper can only be 
  confirmed with the assistance of a computer.  To facilitate verification and 
  further exploration on the reader's part, a commented {\em Mathematica} 
  file on the author's homepage accompanies
  this paper.  It contains some of the formulas and data presented here.    
          
  \subsection{Relation to a similar paper.}
      The polynomials $F_1(p,q,x)$ and $F_2(a,b,x)$ are similar in nature
  to the polynomials $g_{27}(u,v,x)$ and $g_{28}(u,v,x)$ of \cite{CoversM05} which 
  have Galois groups $W(E_6)$ and $W(E_7)^+$ respectively.
  However \cite{CoversM05} and this paper focus on different theoretical topics to avoid 
  duplication.  The discussion of monodromy and the universality of specialization
  sets in \cite{CoversM05} applies after modification to the new base schemes $U_{3,1,1}$
  and $U_{3,2}$ here.  Similarly, our general discussion of
  division polynomials here  could equally well be illustrated by $g_{27}(u,v,x)$ and $g_{28}(u,v,x)$ .
  
  \subsection{Acknowledgements}  It is a pleasure to thank Zhiwei Yun for a conversation about $G_2$-rigidity 
 from which this paper grew.  It is equally a pleasure to thank 
 Michael Dettweiler and Stefan Reiter for helping 
 to make the direct connections to their work \cite{DR}.  We are also grateful
 to the Simons Foundation for research support through grant
 \#209472.

  \section{Background}
\label{background}

    This section provides some context for our later considerations.   

\subsection{Division Polynomials}

     Classical formulas \cite[p.\ 200]{We} let one pass directly from an elliptic curve $X: y^2 = x^3 + a x + b$ 
 to division polynomials giving $x$-coordinates of their $n$-torsion points.  Initializing
 via $f_1 = 1$ and $f_2 = 2$, these division polynomials $f_n$ for $n \geq 3$ are
 computable by recursion:
 \begin{eqnarray*}
 f_{3} & = & 3 x^4 +6 a x^2+12 b x -a^2, \\
 f_4 & = & 4 x^6 + 20 a x^4 + 80 b x^3 -20 a^2 x^2 -16 a b x -4 a^3-32 b^2, \\
 f_{2m} & = & f_m \left( f_{m+2} f_{m-1}^2 - f_{m-2} f_{m+1}^2 \right)/2, \\
 f_{4m+1} & = & (x^3+a x + b)^2 f_{2m+2} f_{2m}^3 - f _{2m-1} f_{2m+1}^3, \\
 f_{4m+3} & = & f_{2m+3} f_{2m+1}^3 - (x^3 + a x + b)^2 f_{2m} f_{2m+2}^3.
 \end{eqnarray*}
 Special cases give interesting number fields.   For example, at $(a,b) = (-1/3,19/108)$ the degree sixty polynomial 
 $f_{11} \in \Q[x]$ has Galois group  $GL_2(11)/\{\pm 1\}$ and field discriminant
  $-11^{109}$. 
 
     On an abstract level, there are interesting number fields
from $n$-torsion points on any abelian variety over $\Q$.  More
 generally, from any variety $Y$ over $\Q$ there are interesting field extensions from
 the natural action of $\Gal(\overline{\Q}/\Q)$ on the cohomology groups
 $H^m(Y(\C),\Z/n\Z)$.   However for most pairs $(X,n)$, there is nothing
 remotely as explicit as the above recursion relations.   In fact, there 
 is presently no way at all to produce explicit division polynomials describing these
 fields.

\subsection{The group $SU_3(3).2 \cong G_2(2)$}
\label{group}

The Atlas \cite{Atlas} provides a wealth of
group-theoretic information about the group $SU_3(3).2 \cong G_2(2)$.  In particular, this group has the form 
$\Gamma.2$, where $\Gamma$ has order $6048 = 2^5 3^3 7$ and is the $12^{\rm th}$ smallest 
non-abelian simple group.

\begin{table}[htb]
\[
\begin{array}{| lr     ll   |}
\hline
\multicolumn{4}{|c|}{\mbox{Classes in $\Gamma$}} \\
C& |C| &  \lambda_{28} & \lambda_{36}   \\
\hline
1A&1&  1^{28} & 1^{36}    \\   
2A&63& 2^{12} 1^{4} &  2^{12} 1^{12}      \\          
3A&56&3^{9}1 &3^{12} \\  
3B&672&3^{9}1&  3^{11}1^{3}  \\      
4AB&2 \cdot 63 &4^{6}1^{4}&4^{6}2^{6}  \\   
4C &378& 4^6 2^2 & 4^6 2^4 1^4    \\   
6A&504& 6^{4}31   &6^{4}3^{4} \\   
7AB&2 \cdot 864 &     7^4 & 7^51     \\
8AB&2 \cdot 756 & 8^{3} 2^{1}1^2&  8^{3}4^{3}  \\      
12AB& 2 \cdot 504 &  12^{2} 31  & 12^{2}6^{2}  \\  
\hline
\end{array}
\;\;\;\;\;\;\;\;\;\;\;\;
\begin{array}{|lrll|}
\hline
\multicolumn{4}{|c|}{\mbox{Classes in $\Gamma.2-\Gamma$}} \\
C & |C| & \lambda_{28} & \lambda_{36}  \\ 
\hline
&&&\\
2b& 252 & 2^{12} 1^{4} &   2^{16} 1^{4}   \\
&&&\\
&&&\\
4d& 252 &4^{6}1^{4}&4^{6}2^{6} \\      
&&&\\
6b&2016& 6^{4} 31 &   6^{5}3^{1}2^{1}1  \\ 
&&&\\
8c&1512&     8^3 4  & 8^3 4^2 2 1^2 \\
12cd&2 \cdot 1008 &  12^{2} 31  & 12^{2}6^{2}   \\
\hline
\end{array}
\]
\caption{\label{wtable} Information about conjugacy classes in $\Gamma.2$}
\end{table}

   Table~\ref{wtable} presents some of the
information that is most important to us.    The left half corresponds to the $14$ conjugacy
classes in $\Gamma$.  The six classes  $1A$, $2A$, $3A$, $3B$, $4C$, and $6A$ are 
rational, while the remaining classes are conjugate in pairs  over
the quadratic fields $\Q(i)$, $\Q(\sqrt{-7})$, $\Q(i)$, and $\Q(i)$ respectively.  
When one considers the full group $\Gamma.2$, these pairs collapse and
one has $16$ conjugacy classes, ten in $\Gamma$ and six in $\Gamma.2-\Gamma$,
with $12c$ and $12d$ two classes conjugate over $\Q(\sqrt{-3})$. 

    Of particular importance to us is that $\Gamma.2$ embeds as a transitive subgroup of  the
alternating group $A_{28}$.  The cycle partition $\lambda_{28}$ associated to a conjugacy
class is given in Table~\ref{wtable}.   The group $\Gamma.2$ also embeds 
as a transitive subgroup of $A_{36}$ and the corresponding $\lambda_{36}$ are given.  
We use the degree $36$ embedding only occasionally.  For example, 
it is useful for distinguishing $3A$ from $3B$ via cycle partitions.   As a convention,
 if we do not refer explicitly to degree we are working
with the degree twenty-eight embedding.

   As just discussed, Table~\ref{wtable} has information about permutation
representations of $\Gamma.2$.  We are also interested in linear representations,
and some group-theoretic information is contained in the small tables at the end of 
\S\ref{UL} (for characteristic 3), at the end of
\S\ref{SL} (for characteristic 2), and in Figure~\ref{g2classpict} 
(for characteristic zero).

\subsection{Rigidity and covers}  Some fundamental aspects of our general 
context are as follows.  
Let $G$ be a finite centerless group and let $C=(C_1,\dots,C_z)$ be a list
of conjugacy classes in $G$.    Define
\begin{eqnarray*}
\overline{\Sigma}(C) & = &  \{(g_1,\dots,g_z) \in C_1 \times \cdots \times C_z : g_1\dots g_z = 1\}, \\
{\Sigma}(C) & = &  \{(g_1,\dots,g_z) \in \overline{\Sigma}(C)  : \langle g_1, \dots, g_z \rangle = G\}. 
\end{eqnarray*}
The group $G$ acts on these sets by simultaneous conjugation and the action is free on ${\Sigma}(C)$.  
The mass of $C$ is $\overline{\mu}(C) := |\overline{\Sigma}(C_1,\dots,C_z)|/|G|$.  
A classical formula, presented in e.g.\ \cite[Theorem~5.8]{MalMat},
gives the mass as a sum over irreducible characters of $G$,
\begin{equation}
\label{massformula}
\overline{\mu}(C) = \frac{|C_1| \cdots |C_z|}{|G|^2} \sum_{\chi} \frac{\chi(C_1) \cdots \chi(C_z)}{\chi(1)^{z-2}}.
\end{equation}
We say that $C$ is rigid if $\mu(C) := |\Sigma(C)|/|G| = 1$ and strictly rigid if moreover $\overline{\mu}(C)=1$.

Let $G \subseteq S_n$ now be a transitive permutation realization of $G$ such that the centralizer of $G$ in $S_n$ 
is trivial.   Let $\tau_1$, \dots, $\tau_z$ be distinct points in the complex projective
line, connected by suitable paths to a fixed basepoint.  
  A tuple $(g_1,\dots,g_z) \in \Sigma(C)$ then
determines
a degree $n$ cover of the projective line with monodromy group $G$ and local monodromy
transformation $g_i$ about the point $\tau_i$.     The genus $g_n$ of the
 degree $n$ cover is calculated via the cycle partitions $\lambda_i \vdash n$
 by the general formula 
 \begin{equation}
 \label{genusformula}
  |\lambda_1|+\dots +|\lambda_z | = (z-2) n+2-2g_n.
  \end{equation}
  Here $|\lambda_i|$  indicates the number of parts of $\lambda_i$.  

    Let $M_{0,w}$ be the moduli space of $w$ labeled distinct points in the projective line.  This
 is a very explicit space, as $\tau_1$, $\tau_2$, and $\tau_3$ can be uniquely normalized
 to $0$, $1$, and $\infty$ respectively.   The group $S_{w}$ acts on $M_{0,w}$ by permuting
 the points.   If $\nu = (\nu_1,\dots,\nu_r)$ sums to $w$ then we let $S_\nu = S_{\nu_1} \times \cdots \times S_{\nu_r}$
 and put $U_\nu = M_{0,w}/S_\nu$.  
 
 When $C = (C_1,\dots,C_z)$ is rigid and the $\tau_i$ move in $M_{0,z}$, 
 all the covers of the projective line fit together into a single cover of
$M_{0,z+1}$.     Moreover, under
simple conditions as exemplified below, this cover is guaranteed to be defined
over $\Q$.   When $z=3$, the space $M_{0,3}$ is just a single point and $M_{0,4}$
 identified with $\bbP^1-\{0,1,\infty\}$, with $\tau_4$ serving as 
 coordinate.  This case has been extensively treated in the literature.
 When $z \geq 4$ the situation is more complicated and 
 a primary purpose of \cite{CoversM05} and the present paper is to give 
 interesting examples.  When some adjacent $C_i$ coincide, 
 the cover descends to a cover of the corresponding quotient $U_\nu$ of $M_{0,z+1}$.

\subsection{The Malle-Matzat cover}  
 Malle and Matzat computed the cover coming from the strictly rigid genus zero triple  $(4d, \, 2b, 12AB)$ belonging
 to the group $\Gamma.2$.   
 This Malle-Matzat cover is similar, but much simpler, than the covers $\pi_1$ and $\pi_2$ that we are about to consider.
 Accordingly we discuss it here as a model, and use it later as well for comparison.  

 Identifying the degree twenty-eight covering curve $X$ with $\bbP_x^1$, the 
cover $\bbP_x^1 \rightarrow \bbP_t^1$ is then given by 
the explicit degree twenty-eight rational function $t=$
\[
\frac{- \! \left(x^6-6 x^5-435 x^4-308 x^3+15x^2+66 x+19\right)^{\! 4}   \! \left(x^4+20 x^3+114 x^2+68 x+13\right)}{2^2 3^9  \left(x^2+4 x+1\right)^{12} (2x+1)}.
\]
The partitions $\lambda_1 = 4^6 1^4$ and $\lambda_3 = 12^2 3 \, 1$
are visible as root multiplicities of the numerator and denominator respectively.  
Rewriting the equation as 
\begin{equation}
\label{mm}
m(t,x) = 0, 
\end{equation}
the partition $\lambda_2 = 2^{12} 1^4$ 
likewise appears as the root multiplicities of $m(1,x)$.

The discriminant of the monic polynomial $m(t,x)$ is 
\begin{equation}
\label{mdisc} D_m(t) = 2^{576} 3^{630}  t^{18} (t-1)^{12}.
\end{equation}
It is a perfect square, in conformity with the fact that $\Gamma.2$ lies in
the alternating group $A_{28}$.  Thus $D_{m}(t)$ is not useful in seeing how the $.2$ enters
Galois-theoretically.    In fact, the order two quotient group corresponds to 
the extension of $\Q(t)$ generated by 
$\sqrt{t(1-t)}$.   

The general theory of three-point covers says that $X \rightarrow \bbP_t$ 
has bad reduction within the primes dividing $|\Gamma.2|$, namely $2$, $3$, 
and $7$.  A particularly interesting feature of $D_m(t)$ is that it reveals
that in fact the Malle-Matzat cover has good reduction at $7$.  In \cite[\S 8]{ABC},
we explained how the Malle-Matzat polynomial is a division
polynomial for a family of varieties with bad reduction only
in $\{2,3\}$, and this connection explains the good reduction at $7$.

\section{Rigid covers of $U_{3,1,1}$ and $U_{3,2}$}  
\label{rigidcovers}
    This section explains how general theory gives the existence 
of our two main covers $\pi_1 : X_1 \rightarrow U_{3,1,1}$ and
$\pi_2 : X_2 \rightarrow U_{3,2}$ and the cubic relation between them.

\subsection{Five strictly rigid quadruples}   
\label{fivecovers}

For a fixed number of ramifying points $z$ and a fixed ambient group $G$, the mass formula \eqref{massformula} lets one find all $C$
with $\overline{\mu}(C) = 1$.    From any explicit tuple $(g_1,\dots,g_z) \in \overline{\Sigma}(C)$,
one gets $\mu(C) = 1$ or $0$ according to whether $\langle g_1,\dots,g_z \rangle$ is
all of $G$ or not.   Carrying out this mechanical procedure for $z=3$ and $G = \Gamma.2$ 
gives several strictly rigid triples, with only the Malle-Matzat triple  
having genus zero.  For
$z=3$ and $G = \Gamma$, one gets yet more rigid triples.  None of these have genus zero
and some of them appear in Table~\ref{sixteen} below.   

Applying this mechanical procedure for $z=4$ yields the
following result:
\begin{Proposition} 
\label{rigidprop} There are no strictly rigid quadruples in $\Gamma.2$.  Up
to reordering and conjugation by the outer involution of $\Gamma$, there
are five strictly rigid quadruples in $\Gamma$:
\[
{\renewcommand{\arraycolsep}{3pt}
\begin{array}{lllrl}
(3A,  & 3A, & 3A, & 4B): & \mbox{(genus $9$),} \\
(4A, & 4A, & 4A, &4B): & \mbox{(genus $9$),} \\
(2A, & 2A, & 3A, & 4A): &  \mbox{(genus $3$),} 
\end{array}
\;\;\;\;\;\;\;\;
\begin{array}{lllrl}
(4A, & 4A, & 4A, & 2A): &  \mbox{(genus $6$),} \\
&\\
(4A,  & 4A, & 3A, & 3A): & \mbox{(genus $9$).} \\
\end{array}
}
\]
Moreover, there are no other rigid quadruples $C$ with 
$\overline{\mu}(C) < 4$.  
   \qed
\end{Proposition}

The list of all quadruples considered in the process of
proving the proposition makes clear that the
 five quadruples presented stand quite apart from
 all the others.    For the case $G=\Gamma.2$, the 
 quadruples $C$ with the smallest $\overline{\mu}(C)$ 
 are $(4d,2b,2A,2A)$, $(4d,2b,3A,2A)$, $(4d,3d,4AB,2A)$,
 $(2b,2b,3A,2A)$, and $(4d,4d,3A,2A)$.  The
 corresponding $(\mu(C),\overline{\mu}(C))$ are 
 $(0,2.750)$, $(3,3.000)$, $(0,3.375)$, $0,3.500)$, and
 $(3,3.\overline{666})$.   For the case of $G=\Gamma$,
 there are fifteen other $C$ with $\overline{\mu}(C) \in [1,2)$;
 all have $\mu(C) = 0$.  Likewise, there are 
 twelve $C$ with $\overline{\mu}(C) \in [2,3)$; 
 four have $\mu(C) = 0$ and eight have 
 $\mu(C)=2$.   Continuing the trend, there are eight $C$ with
 $\overline{\mu}(C) \in [3,4)$; two have $\mu(C)=0$ while
 six have $\mu(C) = 3$.   In particular, as asserted by the proposition, 
 $\mu(C) = 1$ does not 
 otherwise occur in the range $\overline{\mu}(C) <4$;
 we expect that $\mu(C)=1$ does not occur
 either for $\overline{\mu}(C) \geq 4$.

\subsection{The two covers}  In this subsection, we explain how the left-listed quadruples
in Proposition~\ref{rigidprop} all give rise to the same cover $\pi_1 : X_1 \rightarrow U_{3,1,1}$ while the 
right-listed quadruples both give rise to the same cover $\pi_2 : X_2 \rightarrow U_{3,2}$.   
Figure~\ref{twocol} provides a visual overview of our explanation.  

\subsubsection*{The base variety $M_{0,5}$} Let 
\[
M_{0,5} = \Spec \;\Q \!\! \left[s,t,\frac{1}{s(s-1)t(t-1)(s-t)} \right]
\]
be the moduli space of five distinct ordered points in the projective line.  The description on 
the right arises because the five points can be normalized to $0$, $1$, $\infty$, $s$, $t$ by a
unique fractional linear transformation.  

\begin{figure}[htb]
 \[
\begin{array}{ccc}
\includegraphics[width=2.0in]{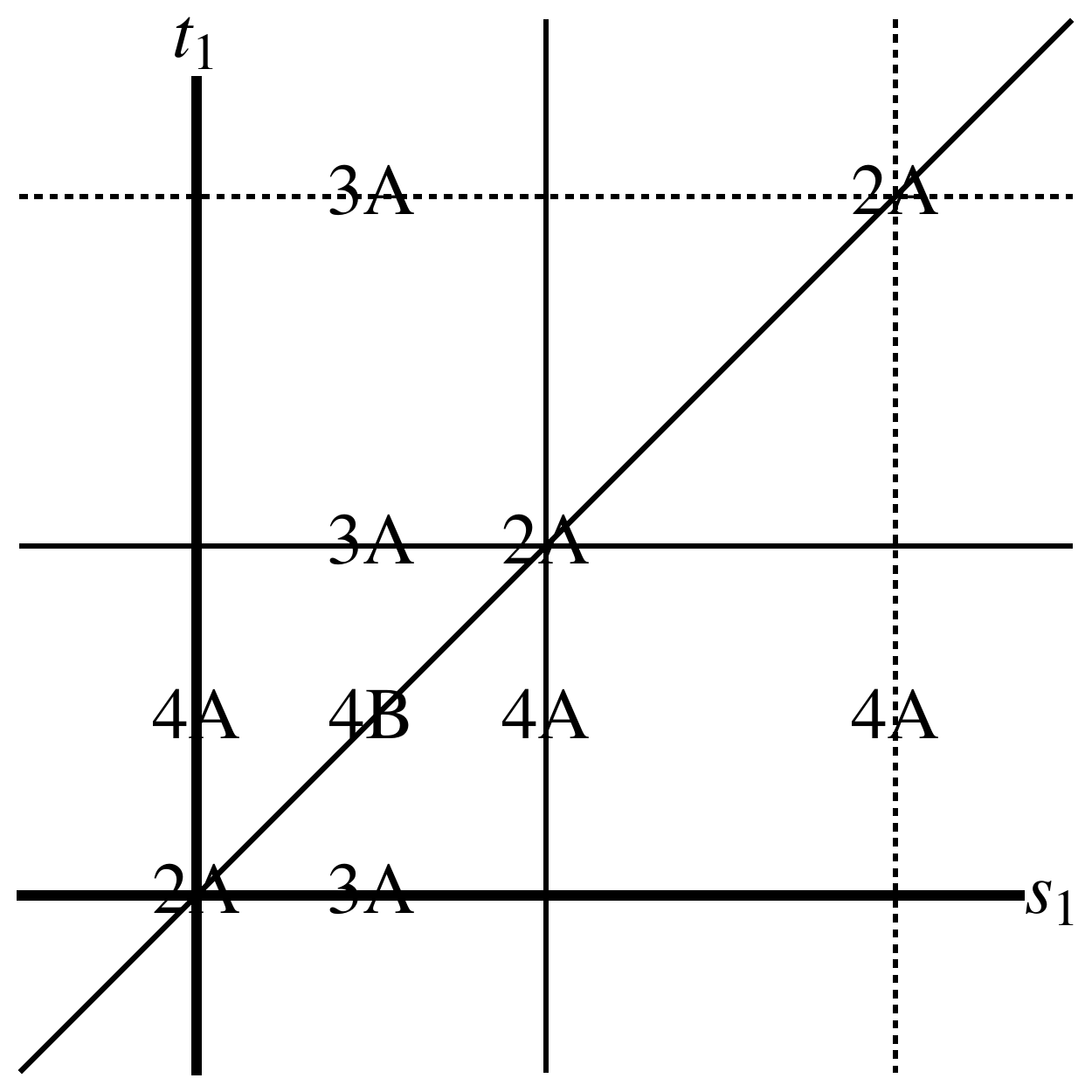} & \;\;\;\;\; \;\;\;\;\;\;\;\; &\includegraphics[width=2.0in]{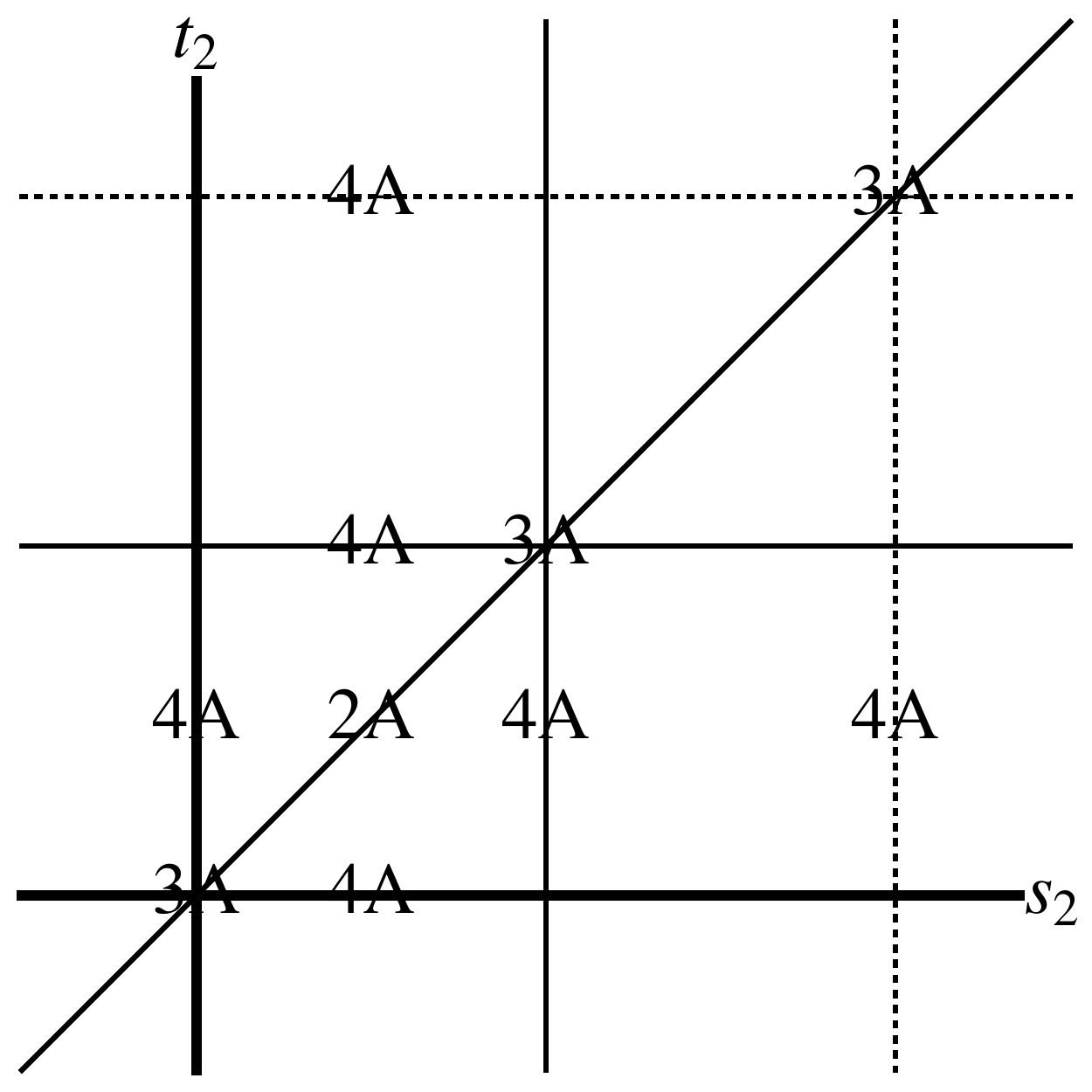} \\
\;\;\;\;\;\;\; \searrow  2 && \; \swarrow 2 \cdot 2 \;\;\;\;\;\;\;  \;\;\;\;\;\;\;\\
\end{array}
\]
\[
\begin{array}{c}
\includegraphics[width=2.0in]{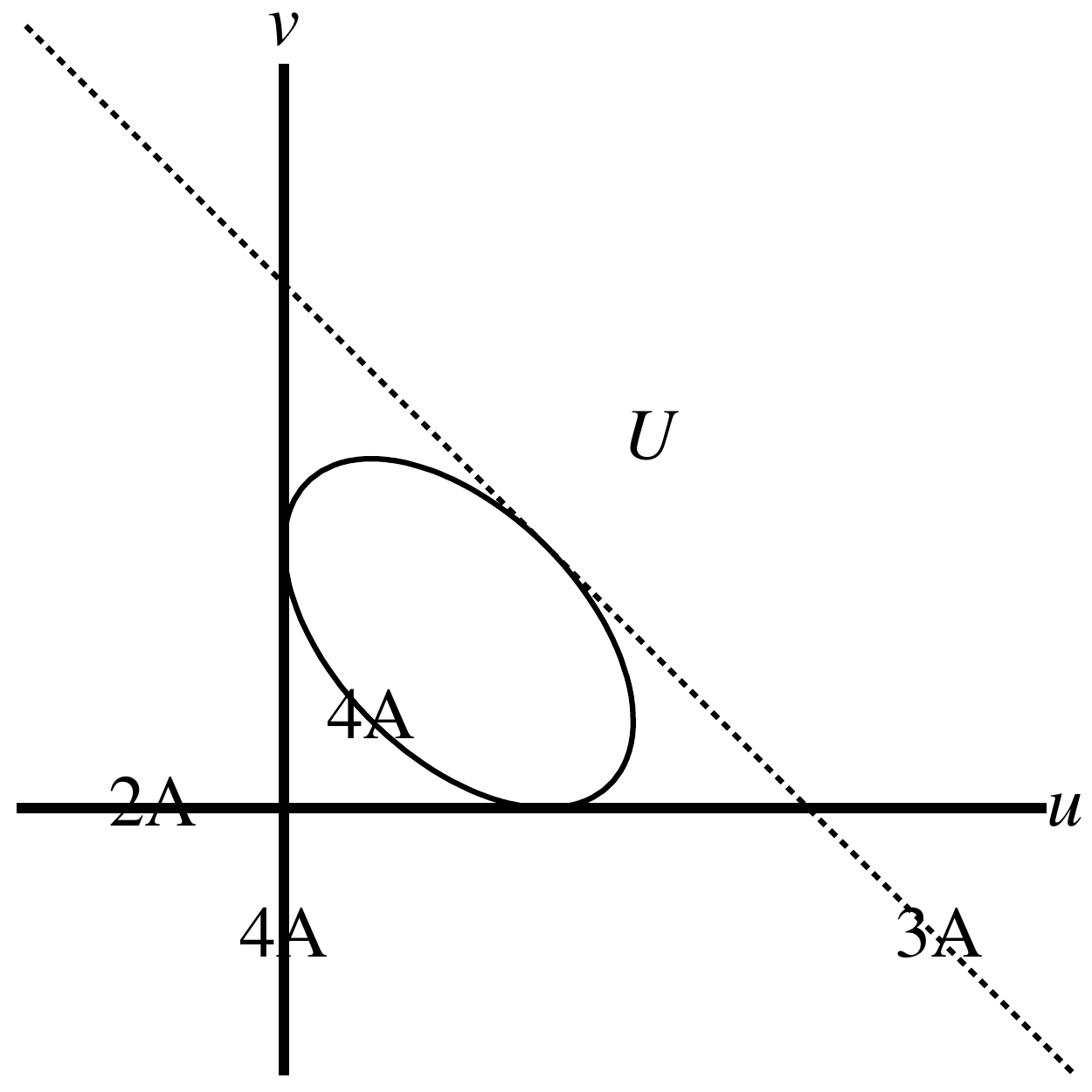}
\end{array}
\]
\[
\begin{array}{ccc}
\;\;\;\;\;\;\;\; \swarrow  3 &  & \; \searrow 3 \;\;\;\;\;\;\; \\
\includegraphics[width=2.0in]{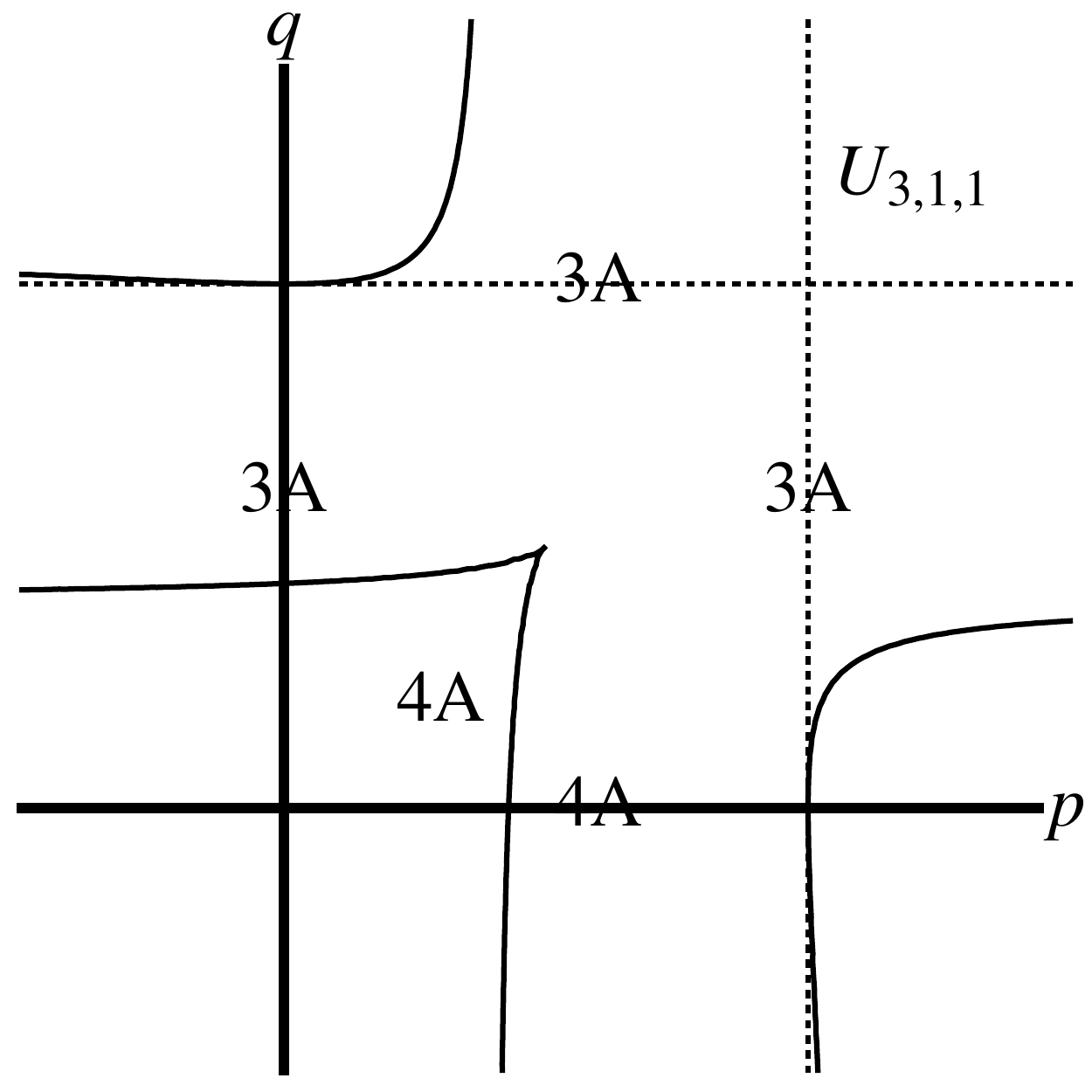} & \;\;\;\;\; \;\;\;\;\;\;\;\; &\includegraphics[width=2.0in]{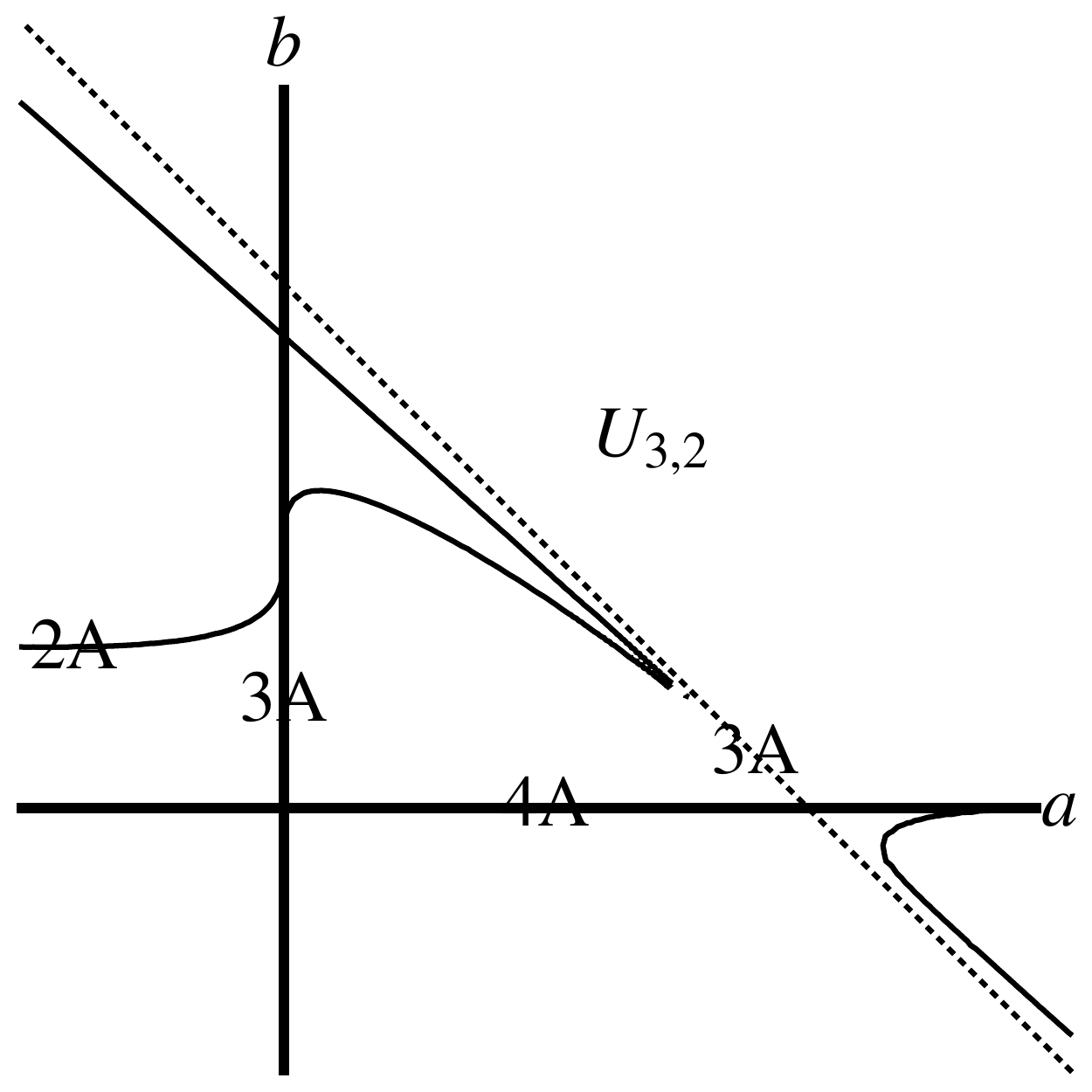} 
\end{array}
 \]
 \caption{\label{twocol}
Base varieties, ramification divisors, and associated conjugacy classes}
 \end{figure}

A naive completion of $M_{0,5}$ is $\overline{M}_{0,5} = \bbP^1_s \times \bbP^1_t$.  The top subfigure 
in each column on Figure~\ref{twocol} gives a schematic representation of the real torus
$\overline{M}_{0,5}(\R)$. As usual, one should imagine the subfigure inscribed in a square,
with the torus obtained by identifying left and right sides, and also top and bottom sides.   
Here and in the rest of Figure~\ref{twocol}, coordinate axes are distinguished by darker lines and
 lines which are at infinity in our particular
 coordinates are indicated by dotting.    
 
A more natural completion $\widehat{M}_{0,5}$ of ${M}_{0,5}$ is obtained from blowing up 
$\overline{M}_{0,5}$ at the three triple points $(0,0)$, $(1,1)$, and $(\infty,\infty)$.  
The natural action of $S_5$ on $M_{0,5}$ extends uniquely 
to $\widehat{M}_{0,5}$.   Reflecting this equivariance, lines in $\widehat{M}_{0,5}-M_{0,5}$ 
are naturally labeled by two-element  
subsets of $\{0,1,\infty,s,t\}$.  Another reflection of equivariance
is that elements of $\{0,1,\infty,s,t\}$ index fibrations
over genus zero curves.    The fibrations $p_s$ and $p_t$
are projections to the $t$ and $s$ axes respectively.  
The smooth fibers of $p_0$ and $p_1$ are the lines going through
$(1,1)$ and $(0,0)$ respectively of slope different from $0$, $1$, $\infty$.  The smooth fibers of $p_\infty$ are certain
hyperbolas going through both $(0,0)$, $(1,1)$.  Note how each fibration partitions the ten lines
of $\widehat{M}_{0,5}-M_{0,5}$ into four sections and
six half-fibers, the latter coming in three pairs to form the singular fibers.

 Consider $(01)$, $(st)$, and $(01\infty)$ in
 their action on $M_{0,5}$.  
 The group $S_3 \times S_2$ that they generate
 acts on the naive compactification $\overline{M}_{0,5}$.
 This action can be readily visualized in terms of our pictures
 of $M_{0,5}(\R)$:  $(01)$ is a half-turn about the point $(1/2,1/2)$,
 $(st)$ is a reflection in the diagonal line, and 
 $(01\infty)$ is a simultaneous one-third turn 
 of the coordinate circles  $\bbP^1_s(\R)$ and $\bbP^1_t(\R)$.

\subsubsection*{Quotients of $M_{0,5}$}  Figure~\ref{twocol}
schematically indicates five planes, each with their 
own coordinates, as indicated by axis-labeling. 
The maps between these planes have the
degrees indicated in Figure~\ref{twocol}, and 
are given by the following formulas:
\begin{align*}
(u,v) & = \left(\frac{(s_1-1) s_1}{(t_1-1) t_1},  \frac{(s_1-t_1)^2}{(t_1-1) t_1} \right),&    (u,v) & =  \left((s_2 - t_2)^2,(s_2 + t_2-1)^2 \right),   \\
&\\
 (p,q) & = \left( \frac{3 (2 u-v+1)}{(u-v+2)^2}, \frac{3 u (u-v+2)}{(2 u-v+1)^2}\right), \!\! &  (a,b) & = \left(\frac{-768 u^3}{W^2}, \frac{9\Delta}{W} \right).
\end{align*}
Here $\Delta = u^2+v^2+1-2 u - 2 v - 2 u v$ is a quantity which
will play a recurring role, while $W = u^2-10 u v+6 u+9 v^2-18 v+9$
is a quantity which appears explicitly here only.    Two moduli interpretations
of $(u,v)$, identifying $U$ with $U_{2,1,1,1}$ and $U_{2,1,2}$ respectively, 
are given in \eqref{x1uv} and \eqref{x2uv} below.   The moduli
interpretation of $(p,q)$ appears in \eqref{x1pq} and \eqref{x2pq} 
below.   The moduli interpretation of $(a,b)$ is less direct, but
arises from the relation \eqref{ab} below.   The four maps
 displayed above are consequences of these
moduli relations.  

Our considerations are mainly birational, and so it 
is not of fundamental importance how we complete
the various planes.  As the diagrams indicate,
three times we complete to a product  
$\bbP^1 \times \bbP^1$ of projective lines, while twice we complete
to a projective plane $\bbP^2$.   We are starting with
two copies of the same variety, with $U^i_{1^5}$ 
having coordinates $s_i$ and $t_i$.   The other varieties 
are quotients:
\begin{align*}
U & = U^1_{1^5}/\langle (01) \rangle, & U & = U^2_{1^5}/\langle (01),(st) \rangle, \\
U_{3,1,1} & = U^1_{1^5}/\langle (01),(01\infty) \rangle, & U_{3,2} & = U^2_{1^5}/\langle (01),(01\infty),(st) \rangle.
\end{align*}
Blowing up some of the 
intersection points would yield more 
natural completions, but we will not be pursuing
our covers at this level of detail.  

The natural double cover $U_{3,1,1} \rightarrow U_{3,2}$ is given in our coordinates by 
\begin{equation}
\label{ab} (a,b) = \left(p^2 q^2 - 6pq +4p+4q-3,p q \right).
\end{equation}
Inserting this map on the bottom row of Figure~\ref{twocol} would of course
make the bottom triangle not commute, as even degrees would be wrong.  
Because of this lack of commutativity, the behavior of $X_1$ over curves
and points in Figure~\ref{u311pict} is not directly related to the behavior of 
$X_2$ over the pushed-forward curves and points in Figure~\ref{pictab}.

\subsubsection*{Covers of $M_{0,5}$.}
The five rigid tuples of Proposition~\ref{rigidprop} enter Figure~\ref{twocol}
 through 
our associating conjugacy classes in $\Gamma$ to lines.  
On the top-left subfigure, from any fixed choice of 
$s \in \C-\{0,1\}$ one has a cover of $\bbP^1_t(\C)$ ramified at $0$, $1$, $\infty$,
and $s$.   The local monodromy classes associated to moving
in a counter-clockwise loop in the $t$-plane about these
singularities form the ordered quadruple $(3A, 3A, 3A, 4B)$.  
On the top-right subfigure they form $(4A, 4A, 4A, 2A)$.  

But now by rigidity one has local monodromy classes associated
to all ten lines of $\widehat{M}_{0,5}-M_{0,5}$.   Using
the monodromy considerations of \cite{CoversM05},
 we have computed these classes. 
The classes
are placed in the top two subfigures of Figure~\ref{twocol}.  Interchanging
the roles of $s$ and $t$, one sees that the cover
of $M_{0,5}$ indicated by the top-left subfigure also
arises from $(4A,4A,4A,4B)$.  However the top right cover now 
just arises in a new way from the original quadruple $(4A,4A,4A,2A)$.    
Via any of the three remaining projections $p_0$, $p_1$ 
$p_\infty$, the covers represented by the top-left and
top-right subfigures arise respectively from 
$(2A,2A,3A,4A)$ and $(4A,4A,3A,3A)$.

\subsubsection*{Descent to covers of $U_{3,1,1}$ and $U_{3,2}$}   The labeling by 
conjugacy classes on both the top-left and top-right 
copies of $M_{0,5}$ is visibly 
stable under the action of 
$S_3 = \langle (01),(01\infty) \rangle $.  Moreover on the top-right,
the labeling is also stable under the diagonal reflection 
$(st)$.    One therefore has descent, to a 
cover $\pi_1 : X_1 \rightarrow U_{3,1,1}$ on the left and 
a cover $\pi_2 : X_2 \rightarrow U_{3,2}$ on the right.

\subsection{Summarizing diagram}
    
\begin{figure}[htb]
\[
\xymatrix@R=1.5pc@C=.8pc{
& \;\; \;X_0 \; \;\;  \ar@{->}[dl]_{\Sigma_1} \ar@{->}[dr]^{\Sigma_2}  \ar@{->}[ddd]_{\pi_0}&  & &      &    \!\!\!\! \Q(x_0,y_0) \!\!\!\! \ar@{-}[dl] \ar@{-}[dr]   \ar@{-}[ddd]&   \\
X_1 \ar@{->}[ddd]_{\pi_1} & & X_2 \ar@{->}[ddd]^{\pi_2} && \!\!\!\!\Q(x_1,y_1)\!\!\!\!  \ar@{-}[ddd] & & \!\!\!\! \Q(x_2,y_2) \!\!\!\!  \ar@{-}[ddd] & \\
&&&&&&\\
&U  \ar@{->}[dl]^{\sigma_1} \ar@{->}[dr]_{\sigma_2} &&& &\!\!\!\! \Q(u,v) \!\!\!\! \ar@{-}[dl] \ar@{-}[dr]  & \\
U_{3,1,1}&& U_{3,2} & \;\;\;\;\;\;\;\; &\!\!\!\! \Q(p,q) \!\!\!\! && \!\!\!\! \Q(a,b) \!\!\!\! \\
}
\]
\caption{\label{coverdiagram1} Left: The covers $\pi_1$ and $\pi_2$, as related by the cover $\pi_0$.  Right: Corresponding function fields.}
\end{figure}
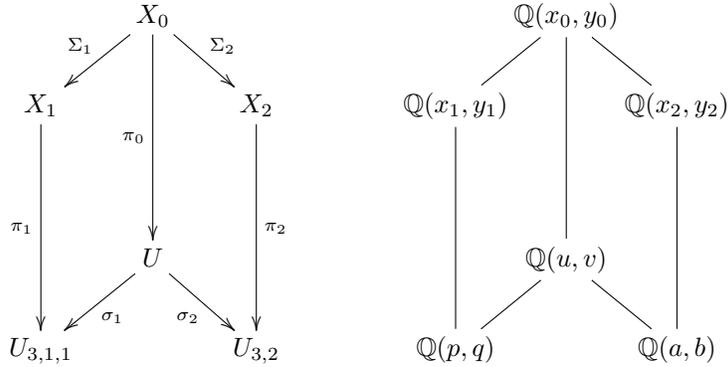

    We now shift attention from Figure~\ref{twocol} to Figure~\ref{coverdiagram1}.  
    The lowest varieties $U_{3,1,1}$, $U_{3,2}$ and their common cubic covering
    by $U$ from Figure~\ref{twocol} are redrawn in the left part of Figure~\ref{coverdiagram1}.  
    The two copies of $M_{0,5}$ from the top of Figure~\ref{twocol} now play a secondary role
    and are suppressed.  In their place, the degree twenty-eight coverings 
    $X_1$ and $X_2$ discussed above are now explicitly indicated.  Also Figure~\ref{coverdiagram1} contains
    their common base change to $X_0 \rightarrow U$.

 The left part of Figure~\ref{coverdiagram1}
    commutes, and so the upper maps $\Sigma_i$, like the lower maps
    $\sigma_i$ from Figure~\ref{twocol}, have degree three.   Note that while
    the top-left part of Figure~\ref{coverdiagram1} has been canonically
    defined, we do not yet have an explicit description for any of the
    surfaces or maps.  We do not yet have an explicit
    description of the vertical maps $\pi_i$ either.   In
    particular, we have not yet discussed the coordinates
    $x_i$, $y_i$ from the top-right part of Figure~\ref{coverdiagram1}.

\section{$3$-division polynomials of Deligne-Mostow covers}
\label{DM} 

    Here we first recognize $\pi_1 : X_1 \rightarrow U_{3,1,1}$ and 
 $\pi_2 : X_2 \rightarrow U_{3,2}$ as associated to $3$-division points
 on certain Deligne-Mostow covers.  Using this connection, we compute $\pi_2$ 
 directly and then deduce explicit formulas for $\pi_0$ and $\pi_1$.  
 The last subsection calculates some sample $L$-polynomials 
 and illustrates how their mod $3$ reductions are 
 determined by our equations for the $\pi_i$.   

\subsection{Local monodromy agreement} Deligne and Mostow's treatises \cite{DM1,DM} concern curves
presented in the form $y^n = f(\mbox{parameters},x)$ and the dependence of their period integrals on the
 parameters.   Their table in \S14.1 of \cite{DM1} has thirty-six lines, each corresponding to a family.  
Their lines $3$ and $2$,  
 written using our parameters $p$ and $q$ on $U_{3,1,1}$, are 
 \begin{eqnarray}
 \label{x1pq} y^4 & = & x^2 (p x^3 + 3 x^2 + 3 x + q), \\
  \label{x2pq} y^4 & = & x  (p x^3 + 3 x^2 + 3 x + q)^2. 
 \end{eqnarray}
 In both cases, the complex roots of $f(x)$ are 
the three roots $\alpha_1$, $\alpha_2$, and $\alpha_3$ of $p x^3 + 3 x^2 + 3 x + q$ and $\alpha_4=0$.   A series solution for 
each equation in the variable $x-\alpha_4 = x$ is 
\begin{align*}
\label{sols}  
y  & =  q^{1/4} x^{1/2} \left(1 + \frac{3x}{4q} - \frac{27x^2}{32 q^2} + \cdots \right), & 
y & =  q^{1/2} x^{1/4} \left( 1 + \frac{3x}{2q} - \frac{9 x^2}{8q^2} + \cdots \right).
\end{align*}
The important quantity for us is the leading exponent associated with $\alpha_4$, namely
$\mu_4 = 1/2$ and $\mu_4 = 1/4$ in the two cases.   Similarly, expanding in the local coordinates
$x-\alpha_i$ for $i \in \{1,2,3\}$, one has $\mu_1=\mu_2=\mu_3$.  In the two cases, these exponents are $1/4$
and $1/2$ respectively.  

Corresponding to the title of \cite{DM} containing just $PU(1,n)$ rather than more general $PU(m,n)$, Deligne
and Mostow are interested in the case when the sum of the $\mu_j$ corresponding to roots of $f(x)$ 
is in $(1,2)$.  
A leading exponent at $\infty$, here $\mu_5$, is then defined so that the sum of all $\mu_i$ is $2$.  
So, summarizing in the two cases, the exponent vector is
\begin{align*}
(\mu_1,\mu_2,\mu_3,\mu_4,\mu_5) & = \left(\frac{1}{4},\frac{1}{4},\frac{1}{4},\frac{1}{2},\frac{3}{4}\right), & 
(\mu_1,\mu_2,\mu_3,\mu_4,\mu_5) & = \left(\frac{1}{2},\frac{1}{2},\frac{1}{2},\frac{1}{4},\frac{1}{4}\right).
\end{align*}
These are the quantities actually presented on lines 3 and 2 of the Deligne-Mosow table.  
From $\mu_4=\mu_5$ one has descent from $U_{3,1,1}$ to $U_{3,2}$ in the second
case, but not the first.  

Switch notation to $(\mu_0,\mu_1,\mu_\infty, \mu_s,\mu_t)$ to agree with the previous section.  The local monodromies
about the divisor of $D_{jk}$, as classes in $GL_3(\C)$, are represented by  
\[
m_{jk} = \left(
\begin{array}{ccc} 
1 & 0 & 0  \\
0 & 1 & 1 \\
0 & 0 & \exp(2 \pi i (\mu_j + \mu_k))
\end{array}
\right).
\]
Here the off-diagonal $1$ can be replaced by $0$, except in the
case $\mu_j+\mu_k \in \Z$, i.e.\ $\mu_j+\mu_k = 1$.  

Global monodromy is in fact in a unitary subgroup of $GL_3(\Z[i])$.  The matrix $m_{jk}$ has infinite order if
$\mu_j+\mu_k=1$, and otherwise has the finite order $\mbox{denom}(\mu_j+\mu_k) \in \{2,4\}$.  
 Reducing
to $PU_3(\F_3) \subset PGL_3(\F_9)$, the infinite order $m_{jk}$ acquire order $3$ and the
finite order $m_{jk}$ maintain their order.  Moreover, not just the orders but 
even the conjugacy classes 
can be shown to agree with those presented at the top of
 Figure~\ref{twocol}.  Thus
the rigid covers of the previous section are realized as 
$3$-division covers.

\subsection{Explicit equations}  The following theorem  gives 
equations describing the three covers $\pi_0$, $\pi_1$, and $\pi_2$. 
A preliminary comment about the contrast between curves and 
surfaces in order.  
Requiring automorphisms to fix $\C$ pointwise,  $\Aut(\C(x))$ is just $PGL_2(\C)$ while $\Aut(\C(x,y))$ is the 
infinite-dimensional Cremona group.  
A consequence is that any given $F : \bbP^1 \rightarrow \bbP^1$
is already in a good form.  Furthermore,  one has only a total of six degrees
of freedom in adjusting domain and target coordinates in
order to get a particularly nice form, like that of the
Malle-Matzat cover.  However a given 
$F: \bbP^2 \rightarrow \bbP^2$ may be in far from best form, 
and adjusting coordinates to improve the form seems to
be more of an art than a science.  The theorem
gives the best form we could find in each case, but does not
exclude the possibility that there
are more concise forms. 

\begin{Theorem}   \label{ExplicitF} The surfaces  $X_0$, $X_1$, and $X_2$ are all rational.  
There are coordinate functions $(x_i,y_i)$ on $X_i$ so that the top
five maps in the left half of Figure~\ref{coverdiagram1} are as follows:
\smallskip

\noindent {\bf The three covers with domain $X_0$.}  Abbreviate $(x,y) = (x_0,y_0)$ and
\begin{eqnarray*}
 g_4 & = & 15 x^2-4 y x-4 x+5, \\
 g_{6a} & = & 9 x y^2+y^2+18 x y-18 y-66 x+6, \\
 g_{6b} & = & 225 x^2-30 y x-30 x-2 y^2+6 y+33, \\
 g_{7a}& = & 15 y x^2+65 x^2-2 y^2 x-4 y x-2 x+5 y+5, \\
 g_{7b}& = & 45 y x^2-105 x^2+6 y^2 x-8 y x-14 x-5 y-5, \\
 g_9& = & 225 x^3-30 y x^2-105 x^2+y^2 x+22 y x+21 x+y^2-8 y-9, \\
 g_{10} & = & 225 x^2 y^2+1200 x^2 y+2850 x^2+250 x y^2-1500 x+2 y^4+8 y^3+37 y^2 \\&& \qquad -192 y+402, \\
 g_{17} & = & 2025 x^3 y^2+24300 x^3 y+39150 x^3+540 x^2 y^3+1845 x^2 y^2-180 x^2 y \\ && \qquad  -29610 x^2-18 x
    y^4+168 x y^3-213 x y^2-252 x y+9522 x+10 y^4 \\ && \qquad +60 y^3-105 y^2+900 y-2070, \\
 g_{18} & = & 50625 x^5-3375 y x^4+30375 x^4-675 y^2 x^3+2025 y x^3+2700
      x^3+75 y^3 x^2 \\ && \qquad +2025 y x^2+5850 x^2-2 y^4 x-18 y^3 x+33 y^2 x-513 y x+63
      x \\ && \qquad +15 y^3 -30 y^2+270 y+315. 
\end{eqnarray*}
Then $\Sigma_1$, $\pi_0$, and $\Sigma_2$ are given by 
\begin{align*}
x_1 & =  -\frac{g_{17}}{g_{6a} g_{6b}}, & y_1 & = \frac{45 g_{10}  (x+1) \left(9
    x^2-2 x+1\right)}{g_4 g_{6a}  g_{6b}}, \\
  u & =  
  \frac{g_{6b}^6 x^4 (x+1)}{25 g_{6a}^3 g_{7a}^2 \left(9 x^2-2   x+1\right)}, & v & = 
   \frac{g_{18}^2 g_4^2 g_9}{25 g_{6a}^3 g_{7a}^2 \left(9    x^2-2 x+1\right)}, \\
x_2 & =  \frac{1}{x+1},  & y_2 & = 
 \frac{g_4  g_{7a}}{5 g_{7b} (x+1)^2}.
\end{align*}
\noindent {\bf The cover $\pi_1$.}   Abbreviate $(x,y) = (x_1,y_1)$ and  
\begin{eqnarray*}
h_{11} & = & 2 x^5+4 x^4-12 x^3 y+8 x^3+9 x^2 y+16 x^2-24 x y+8 x+3 y^3+18 y+16, \\
h_{26} & = & 4 x^7+48 x^6 y+7 x^6-72 x^5 y^2+24 x^5+48 x^4 y^3-63 x^4 y^2+288 x^4 y \\ && +42 x^4-12 x^3
    y^4-288 x^3 y^2+48 x^3-3 x^2 y^4+192 x^2 y^3-252 x^2 y^2 \\ && +576 x^2 y+84 x^2-24 x y^4-288 x y^2+32 x+3
    y^6-6 y^4+192 y^3-252 y^2 \\ && +384 y+56.
\end{eqnarray*}
Then $\pi_1$ is given by 
\begin{align*}
p & =  \frac{3 h_{26}}{h_{11}^2}, &
q & =  \frac{3 h_{11} y \left(3 x^2-y^2+6\right)^4}{h_{26}^2}.
\end{align*}
\noindent {\bf The cover $\pi_2$.}  Abbreviate $(x,y) = (x_2,y_2)$ and
\begin{eqnarray*}
f_9 & = & 144 x^3 y-408 x^2 y-12 x^2+8 x y^2+388 x y+20 x-9 y^2-126 y-9, \\
f_{14} & = & 36 x^4 y^2-288 x^3 y^2-504 x^3 y+816 x^2 y^2+1236 x^2 y-12 x^2+2 x y^3 \\ && -840 x y^2-1038
    x y+20 x-9 y^3+297 y^2+297 y-9.
\end{eqnarray*}
Then $\pi_2$ is given by 
\begin{align*} a & = \frac{3 f_9^3}{f_{14}^2 \left(12 x^2-20 x+9\right)}, &
b & = \frac{36 x^4
    y^2}{f_{14}}.
\end{align*}
\end{Theorem} 

\proof We will sketch our construction only, as there were many complicated variable changes to reduce to the relatively concise
formulation given in the theorem.  We first found $\pi_2$ as follows.  Via \eqref{x2pq}, the genus three curve $Y_2(p,q)$
is presented as a quartic cover of $\bbP^1_x$.     Replacing $x$ by $t^2$ in \eqref{x2pq} and 
factoring, one gets a presentation of $Y_2(p,q)$ as a double cover of the $t$-line $\bbP^1_t$:
\[
y^2 = p t^7 + 3 t^5 + 3 t^3 + q t.
\]
Three-torsion points on the Jacobian of $Y_2(p,q)$ are related to unramified abelian triple covers of $Y_2(p,q)$.  
Such a triple cover arises as a base-change of certain ramified non-abelian triple covers of
$\bbP_t^1$.  

Consider now a partially-specified triple cover of $\bbP^1_t$, given by 
\[
(a t+\mu ) z^3 +(b t- \lambda \mu ) z^2 +t (c+t) z+t (d-t) = 0.
\]
The discriminant of this polynomial with respect to $z$ is a septic
polynomial in $t$ with zero constant term.  Setting it equal to 
$k (p t^7 + 3 t^5 + 3 t^3 + q t)$ imposes the necessary
ramification condition.  It also gives seven equations in the
seven unknowns $a$, $b$, $c$, $d$, $\lambda$, $\mu$, and $k$, 
all dependent on the two parameters $p$ and $q$. 

The equations corresponding to the coefficients of $t^2$, $t^4$, and $t^6$ let 
one eliminate $d$ and $\mu$ and reduce the remaining equation to
\begin{eqnarray*}
\!\!\!\!&\!\!\!\!&-9 a^2 \lambda ^2+162 a^2 \lambda -729 a^2+18 a b \lambda ^2+180 a b \lambda -486 a b+120 a c
    \lambda \\ \!\!\!\!&\!\!\!\!& \qquad \qquad \qquad  -216 a c+3 b^2 \lambda ^2+34 b^2 \lambda +27 b^2+24 b c \lambda +72 b c+32 c^2=0.
\end{eqnarray*}
The system consisting of this equation and the equations coming from the
 coefficients of $t^1$, $t^3$, $t^5$, and $t^7$ is very complicated to solve. 
Nonetheless,  one can eliminate all the remaining variables, at the expense
of putting in the new parameters $x_2$ and $y_2$.   Conveniently, $p$ and $q$ 
enter the final formulas symmetrically and can then be replaced by 
$a$ and $b$ via \eqref{ab}, yielding our presentation of $\pi_2$.

Our formula for $\pi_0$ was then obtained via base-change.  To get $\pi_1$ we first
built a double cover $\tilde{X}_0$ of $X_0$ which is a Galois sextic cover of the
yet-to-be explicitized $X_1$.  Then we explicitized $X_1$ by taking invariants under
the Galois action,  $X_1 = \tilde{X}_0/S_3$.   \qed

In general, suppose given a cover of rational surfaces via say $u = A(x,y)/B(x,y)$ and $v = C(x,y)/D(x,y)$. 
Assuming $x$ indeed generates the field extension, one can express the cover in terms of 
$x$ alone via a 
resultant
\[
F(u,v,x) = \mbox{Res}_y(A(x,y) - u B(x,y),C(x,y) - v D(x,y)).
\]
Carrying this out in our context gives $F_0(u,v,x)$, $F_1(p,q,x)$ and 
$F_2(a,b,x)$.  Expanded out, they have 1606, 772, and 209 terms respectively.
Interchanging the roles of $x$ and $y$, one gets polynomials $G_0(u,v,y)$, $G_1(p,q,y)$,
and $G_2(a,b,y)$ with 4941, 1469, and 951 terms respectively.  Again the setting of
surfaces is much more complicated than that of curves. In general, keeping
either just $x$ or just $y$ is unlikely to minimize the number of terms.
More likely the minimum can only be obtained by keeping some third variable
$z \in \Q(x,y)$. There do not seem to be
standard procedures to find these best variables. 

\subsection{$L$-polynomials of Deligne-Mostow covers and their reduction modulo $3$}
\label{UL}
To explicitly illustrate the $3$-division nature of the main polynomials $F_1(p,q,x)$ 
and $F_2(a,b,x)$, we pursue the polynomial $F_0(u,v,x)$ describing their
common base-change.  
Cubically base-changed to the $u$-$v$ plane, the Deligne-Mostow covers in question after some twisting
 become as follows: 
\begin{align}
\label{x1uv}  Y_1(u,v): && v y^4 & =  x^2 (x-1)^3 (v x^2 + (1-u-v) x + u) && \mbox{(genus four)},  \\
\label{x2uv} Y_2(u,v): &&  4 y^4 & =  (x^2 + 2 x + 1-\frac{4}{v})^2 (x^2-2x+1-\frac{4u}{v}) && \mbox{(genus three)}, \\
\nonumber E(u,v): && y^2 & =  (x-1) (v x^2 + (1-u-v) x + u) && \mbox{(genus one)}.
\end{align}
The quadratic subcover of $Y_1(u,v)$ is the elliptic curve $E(u,v)$ while
the quadratic subcover of $Y_2(u,v)$ has genus zero.  

Our monodromy considerations give a relation between 
$Y_1(u,v)$ and $Y_2(u,v)$.   The twisting factors $v$ and $4$ in the equations
above are included so that we can give a clean statement of 
this relation on a more refined level:
\begin{equation}
\label{Lfact}
L_p(Y_1(u,v),x) = L_p(Y_2(u,v),x) L_p(E(u,v),x).
\end{equation}
Here $u$ and $v$ are rational numbers and
 $p$ is any prime good for all three curves.  Each $L$-polynomial $L_p(Y,x)$ is the numerator of the corresponding zeta-function $\zeta_p(Y,x)$, obtained
by determining the point counts $|Y(\F_{p^f})|$ for $f$ up through $\mbox{genus}(Y)$.   Our computations below obtain this $L$-polynomial
via {\em Magma}'s command \verb@ZetaFunction@ \cite{Mag}.  

The factorization \eqref{Lfact} has the
following explicit form:
\begin{eqnarray*}
L_p(Y_1(u,v),x) & = & (1+a x + b x^2 + c x^3 + p b x^4 + p^2 a x^5 + p^3) (1 + d x + p x^2).
\end{eqnarray*}
For $p \equiv 1 \; (4)$, both factors in turn split over $\Q(i)$ as the product of two conjugate polynomials.
For $p \equiv 3 \; (4)$, the coefficients $a$, $c$, and $d$ all vanish, so that each factor is 
an even polynomial.  Taking $(u,v) = (-4,-3)$ as a running example, these two 
cases are represented by the first two good primes:
\begin{eqnarray*}
L_5(Y_1(-4,-3),x) & = &  \left(1 - x^2 - 16 x^3 - 5 x^4 + 125 x^6 \right) \left(1-2 x + 5 x^2 \right), \\
                              & = & N \! \left(1 + ix - (1 - 2i) x^2 - (10 + 5i) x^3 \right) N \! \left(1-(1+2i)x   \right), \\
L_7(Y_1(-4,-3),x) & = & \left(1 + 5 x^2 + 35 x^4 + 343 x^6\right)  \left(1+7 x^2\right).
\end{eqnarray*}
Here and below, $N(f) = f \bar{f}$ is the product of a polynomial $f$ and its conjugate $\bar{f}$.

Consider now $L_p(Y_2(u,v),x) = N(1 + \alpha x + \beta x^2 + \gamma x^3)$ in $\F_3[x]$ 
for varying $p \equiv 1 \; (4)$.   To twist into a situation governed by $SU_3(3)$ we
replace $x$ by $-\gamma^5 x$ to obtain the modified polynomial 
$\hat{L}_p(Y_2(u,v),x) = N(1- \alpha \gamma^5 x + \beta \gamma^2 x - x^3) \in \F_3[x]$.  
Similarly consider $L_p(Y_2(u,v),x) = 1 + b x^2 + b p x^4 + p^3 x^3$ in $\F_3[x]$.
To twist into a situation governed by $SU_3(3).2-SU_3(3)$, we replace
$x^2$ by $p x^2$ obtaining $\hat{L}_p(Y_2(u,v),x) = 1 + b p x^2 + b x^4 + x^6 \in \F_3[x]$.  
For $5 \leq p \leq 97$, the polynomials $\hat{L}_p(Y_2(-4,-3),x)$ are 
calculated directly by \verb@ZetaFunction@ to be
\[
\begin{array}{llll}
\mbox{Class$(p)$} & \lambda_{28}(p) &  \hat{L}_p(Y_2(-4,-3),x) \in \F_3[x] & \mbox{Primes $p$} \\
\hline
3B &3^9 1&  N \! \left(1-x^3\right) & 89\\
7AB & 7^4 & N \! \left(1-(1+i)x + (1-i) x^2-x^3\right) &
   5,13,29,53,61,73,97 \\
8AB & 8^3 2 \, 1^2 & N \! \left(1-ix -ix^2-x^3\right) & 37,41 \\
12AB &  12^2 3 \, 1 & N \! \left(1+(1-i) x -(1+i) x^2-x^3\right)  & 17 \\
6b& 6^4 3 \, 1 &  \left(1+x^2\right)^3 & 11,19 \\
8c &8^3 4 &  \left(1+x^2\right) \left(1+2x+2x^2\right) \left(1+x+2x^2\right) &
   43,67,79,83 \\
   12c, 12d & 12^2 3 \, 1 & (1+x)^2 (1+2x)^2 \left(1+x^2\right) & 7,23,31,47,59,71. \\
\end{array}
\]
For general $(u,v)$, the fact that $F_0(u,v,x)$ functions as a 3-division polynomial 
is seen by the fact that $\hat{L}_p(Y_2(u,v),x) \in \F_3[x]$ depends only 
the conjugacy class in $\Gamma.2$ determined by $p$.
Up to small ambiguities, as described in Table~\ref{wtable}, this
conjugacy class is determined by the class
of $p$ modulo $4$ and the factorization partition
$\lambda_{28}(p)$ of $F_0(u,v,x) \in \F_p[x]$.

\section{$2$-division polynomials of Shioda quartics}
\label{shioda} 
    In this section, we recognize $\pi_1 : X_1 \rightarrow U_{3,1,1}$ and $\pi_2 : X_2 \rightarrow U_{3,2}$ 
 as $2$-division polynomials for certain genus three Shioda curves.  
 %
 The last subsection calculates some sample $L$-polynomials 
 and illustrates how their mod $2$ reductions are 
 determined by our equations for the $\pi_i$.

\subsection{The Shioda $W(E_7)^+$ polynomial}  In \cite{Shi}, Shioda exhibits multiparameter
 polynomials for the Weyl groups $W(E_6)$, $W(E_7)$, and $W(E_8)$.   He proves 
 in Theorem~7.2 that these polynomials
 are generic, in the sense that any $W(E_n)$ extension of a characteristic
 zero field $F$ is given
 by some specialization of the parameters.

The case of $W(E_7) \cong W(E_7)^+ \times C_2$ is explained in greater detail in \cite{Shi2} and goes as follows.  
  Fix a parameter vector $r=(r_1,r_3,r_4,r_5,r_6,r_7,r_9) \in \C^7$ and
consider the 
equation
\begin{equation}
\label{shiodaell}
y^2 = x^3 + (w^3 + r_4 w + r_6) x + (r_1 w^4 + r_3 w^3 + r_5 w^2 + r_7 w + r_9).
\end{equation}
The vanishing of the right side defines a quartic curve $Q_r$ in the
$w$-$x$ plane.  The equation itself defines a $K3$ surface in $x$-$y$-$w$ space 
mapping to the $w$-line
with elliptic curves as fibers.   Now consider the substitutions
\begin{align*}
x & =  z w + b, &
y & =  c w + d w + e,
\end{align*}
which make each side of \eqref{shiodaell} a quartic polynomial in $w$.   Equating
like coefficients, \eqref{shiodaell} then becomes five equations in the five unknowns
$z$, $b$, $c$, $d$, and $e$.   There are fifty-six solutions, paired according
to the negation operator $(z,b,c,d,e) \mapsto (z,b,-c,-d,-e)$.    Much of the
interest in Shioda's theory comes from regarding these solutions as generators
for the rank seven Mordell-Weil group of the generic fiber.   

    Our interest instead is that the twenty-eight lines $x = z w + b$ are exactly the
twenty-eight bitangents of $Q_r$.   
The variables $b$, $c$, $d$, and $e$ can be very easily eliminated and
one gets Shioda's degree twenty-eight generic polynomial for the 
rotation subgroup $W(E_7)^+$:
\begin{eqnarray*}
\lefteqn{S(r,z) = } \\
 & &  z^{28} -8 r_1 z^{27} + 72 r_3 z^{25} +60 r_4 z^{24} + (-504 r_5 + 432 r_1 r_4)  z^{23} + \\
 && (384 r_1^2 r_4-1248r_1 r_5+540 r_3^2-540 r_6) z^{22} +  \cdots
\end{eqnarray*}
Expanded out as an element of $\Z[r,z] := \Z[r_1,r_3,r_4,r_5,r_6,r_7,r_9,z]$, there are 1784 terms. 
 The polynomial is weighted homogeneous when
the variable $z$ is given weight $1$ and each parameter $r_i$ is given weight
$i$.  The polynomial discriminant of $S(r,z)$ factors over $\Q$ as
$
\Delta(r) = D(r) C(r)^2,
$
with $D(r)$ a source of ramification and $C(r)$ an
irrelevant artifact of our coordinates.    

\subsection{Using $\Gamma.2 \subset W(E_7)^+$}
The group $\Gamma.2$ is a subgroup of $W(E_7)^+$.  Since genericity implies descent-genericity \cite{Kemp},
any degree $28$ extension 
$K/F$ with Galois group $\Gamma.2$ is of the form $F[x]/S(r,z)$ for suitable
$r \in F^7$.   For the Malle-Matzat polynomial $m(t,z)$, we considered various $t \in \Q$ 
and conducted a very modest  search over different polynomials of small height
 defining the same field as $\Q[x]/m(t,z)$.  
For a few $t$, we found a polynomial of the form $S(r,z)$ for certain
$r \in \Q^7$.  Some of these seven-tuples had similar shapes, and interpolating
these only we found that the Malle-Matzat family seemed also to be given by
\begin{equation}
\label{mms}
S(0, -27 t^2, -81 t^2, 243 t^3, 243 t^3, -729 t^4, 729 t^5,z)=0.
\end{equation}
The correctness of this alternative equation is algebraically
confirmed by eliminating $t$ from the pair of 
equations \eqref{mm}, \eqref{mms}, to obtain
the relation
\begin{equation}
\label{mmxz}
   z =  \frac{\begin{array}{c}(x-1)^2 \left(x^4+20 x^3+114 x^2+68 x+13\right) \cdot  \\ \left(x^6-6 x^5-435 x^4-308
    x^3+15 x^2+66 x+19\right)^2 \end{array}}{243 \left(x^2+4 x+1\right)^8}.
\end{equation}
Thus Equation~\eqref{mms} realizes the Malle-Matzat polynomial as a $2$-division polynomial
for an explicit  family of genus three curves.  

The simplicity of the equational form \eqref{mms} is striking, especially taking into account that
all the positive integers printed are powers of $3$.   Expanding the family
out as a polynomial in $\Z[t,z]$ hides the simplicity, as
there are $75$ terms. 

\subsection{A search for $\Gamma.2$ specializations}
Given the simplicity of \eqref{mms}, we searched for similar families 
as follows.  We considered one-parameter polynomials of the
form $S(r,z)$ with $r_i = a_i t^{e_i}$.  Here the $e_i \in \Z_{\geq 0}$ are
fixed and the constants $a_i$ yet unspecified.  We looked at many $(e_1,e_3,e_4,e_5,e_6,e_7, e_9)$ 
near-proportional to $(1,3,4,5,6,7,9)$ so as to ensure that 
$D(a_1 t^{e_1},\dots,a_9 t^{e_9})$ has the form $t^a d(t)$ with
$d(t)$ of small degree.   When a particular exponent $e_i$ 
made a proportionality $(e_1,\dots, e_9) \propto (1,\dots,9)$
not so close, we set $a_i$ equal to zero, rendering $e_i$ irrelevant.

We then worked modulo $5$, letting $(a_1,\dots,a_9)$ 
run over relevant possibilities in $\F_5^7$.   If $k$ of the $a_i$ are  
set equal to zero, we  looked at 
 just $4^{5-k}$ possibilities: we keep the other $a_i$ nonzero, and 
homogeneity and the scaling $t \mapsto ut$ each save a
 factor of $4$.  We examined each one-parameter
family $S(a_1 t^{e_1},\dots,a_9 t^{e_9},x)$ 
by specializing to $t \in \F_{5^j}$ and factoring
in $\F_{5^j}[x]$.   In the rare cases when all factorization 
patterns $\lambda_{28}$ for $j=1$, $2$, and $3$
correspond to elements of $\Gamma.2$, as on
Table~\ref{wtable}, we proceeded under the expectation
 that $S(a_1 t^{e_1},\dots,a_9 t^{e_9},x)=0$
defines a cover with Galois group in $\Gamma.2$.  

For fifteen $(e_1,\dots,e_9)$ we found exactly 
one $(a_1,\dots,a_9)$ which works.  For five $(e_1,\dots,e_9)$
we found several $(a_1,\dots,a_9)$ which work, suggestive
of a two-parameter family.    We then reinspected
 these five $(e_1,\dots,e_9)$ in characteristic seven,
imposing also that the covers sought be tame.  
The case $(e_1,e_3,e_4,e_5,e_6,e_7,e_9) = (0,1,1,\star,2,2,2)$
seemed to give a two-parameter family in both 
characteristics, satisfying the tameness condition at $7$;
here the $\star$ means that we are setting $a_5=0$.   
Standardizing coordinates, the two-parameter 
families seemed to match well, and there remained the
task of lifting to characteristic zero.

       We first found that $S(1,0,3 t,0,0,0,-t^2,z) \in \Q[v,z]$ defines
a $3$-point cover, giving us hope that coefficients
might be even simpler than in \eqref{mms}.  
Finally we found a good two-parameter
family $S_0(u,v,z)=0$ where 
\begin{equation}
\label{fdef}
S_0(u,v,z) =  S(1,u-v+1,-3 u,0,u (-u+v-1),u (-u+v-1),-u^2,z).
\end{equation}
The discriminant of $S_0(u,v,z)$ is
\[
D(u,v) = 2^{216} 3^{108} u^{42} v^{24} (u^2-2 u v-2 u+v^2-2 v+1)^2
\]
times the square of a large-degree irreducible polynomial in $\Z[u,v]$.

\subsection{Explicit polynomials}    Our computation of $S_0(u,v,z)$, as just described,
is completely independent of the considerations of the previous section.  In fact
we found $S_0(u,v,z)$ before we found its analog $F_0(u,v,x)$ 
from the previous section.   It might have been possible to directly desecnd
$S_0(u,v,z)$ to $S_1(p,q,z)$ and $S_2(a,b,z)$ below.   However
instead we obtained these new $S_i$ from the corresponding
$F_i$:  we took lots of specialization points, applied {\em Pari}'s \verb@polred@
to obtain alternate polynomials, selected those that are of the form 
$S(r_1,r_3,r_4,r_5,r_6,r_7,r_9,z)$, and interpolated those 
that seemed to fit a common pattern.  
 %

\begin{Theorem}  
\label{ExplicitS} Abbreviate $d = p^2 q^2-6 p q+4 p+4 q-3$, $A=256/a$, and $B = (b-1)/8$.     
The covers $\pi_0$, $\pi_1$, and $\pi_2$
are also given respectively via the polynomials $S_0(u,v,z)$, 
\begin{eqnarray*}
S_1(p,q,z) & = & S(  \begin{array}{r}
                   0, \\
                   d^2 p, \\
                   3 d^2 p^2 (q-1), \\
                   3 d^3 p^2, \\
                   -d^3 p^2 \left(3 p^2 q^2-9 p q+4 q+2 p\right), \\
                   -3 d^4 p^3 (q-1), \\
                   d^5 p^4 \left(2 p q^2-3 q+1\right), 
                   \end{array} z ), \\
 S_2(a,b,z) & = & S(\begin{array}{r}
  1, \\
                   3 \left(A B^2+2\right), \\
                   -3 \left(8 A B^2+A B+1\right), \\
                   -3 \left(5 A B^2+A B-4\right), \\
                   -8 A^2 B^4-A^2 B^3-184 A B^2-31 A B-A-2, \\
                   -56 A^2 B^4-7 A^2 B^3-199 A B^2-58 A B-4 A+10, \\
   \!\!\!\!\!\!\!\!\!\!\!\!\!\!\!\!\!\!\!                -440 A^2 B^4-103 A^2 B^3-6 A^2 B^2-693 A B^2-183 A B-12 A+3,
 \end{array} z).
\end{eqnarray*}
\end{Theorem}

\proof  We describe Case 0, as the other cases are similar except that 
the analog of \eqref{zans} is much more complicated.   Analogously to \eqref{mmxz}, 
One needs to find $z$ in the function field $\Q(x,y)$ of $X_{0}$
satisfying $S_{0}(u,v,z)=0$.  To find a candidate $z$, one 
takes a sufficiently large collection of  $\{(x_i,y_i)\}$ 
of ordered pairs in $\Q^2$.  One next obtains the pairs 
$(u_i,v_i) = \pi_{0}(x_i,y_i)$.  Discarding the very rare cases
where $S_0(u_i,v_i,z) \in \Q[z]$ has more than one rational root, 
one defines $z_i$ to be the unique rational root of $S_0(u_i,v_i,z)$.  
The desired $z$ is then obtained by interpolation, being
\begin{equation}
\label{zans}
z = \frac{(3x-1)f_6}{g_6} = \frac{(3x-1)\left(9 x y^2+18 x y-66 x+y^2-18 y+6\right)}{225 x^2-30 x y-30 x-2
    y^2+6 y+33}.
\end{equation}
Correctness is confirmed by algebraically checking that $S_{0}(u(x,y),v(x,y),z(x,y))$ indeed
simplifies to zero in 
$\Q(x,y)$.  \qed

Fully expanded out, $S_0(u,v,z)$, $S_1(p,q,z)$, and $S_2(a,b,z)$ respectively have
551, 7299, and 1053 terms.   Thus given Shioda's master polynomial $S$, our
$S_i$ admit the relatively concise presentations given in \eqref{fdef} and Theorem~\ref{ExplicitS}.  
Without $S$, the new $S_i$ are of comparable complexity to the previous $F_i$, in
the sense of number of terms.   

\subsection{$L$-polynomials of Shioda quartics and their reduction modulo $2$}
\label{SL}
To illustrate the $2$-division nature of the polynomials $S_0(u,v,z)$, 
$S_1(p,q,z)$, and $S_2(a,b,z)$, one could take any parameter 
pair for which the corresponding polynomial is separable.  As in
\S\ref{UL}, we work with $(u,v) = (-4,-3)$.  

The images of $(u,v)$ in the lower planes are $(p,q) = \sigma_1(-4,-3) = (-12,-3/4)$
and $(a,b) = \sigma_2(-4,-3) = (192,9)$.  By plugging into the three parts of Theorem~\ref{ExplicitS}, and scaling by
$r_i \mapsto r_i/9^i$ in the middle case, one 
gets indices
\[
{\renewcommand{\arraycolsep}{2pt}
\begin{array}{rcrrrrrrr} 
I_0(-4,-3) & \,= \,&  (1, & 0, & 12, & 0, & 0, & 0, & -16), \\
I_1(-12,-3/4) & \, = \,&  (0, & -12, & -84, & -144, & 720, & -1008, & 7872), \\
I_2(192,9) & \,= \, &  (1, & 10, & -39, & -12, & -306, & -450, & -2157). 
\end{array}
}
\]
Taking these vectors as $(r_1,r_3,r_4,r_5,r_6,r_7,r_9)$ and substituting
into the right side of \eqref{shiodaell}, one gets three quartic plane curves,
to be denoted here simply $Q_0$, $Q_1$, and $Q_2$.  

As in \S\ref{UL}, each of the genus three curves $Q_i$ has good $L$-polynomials
\[
L_p(Q_i,x)= 1 + a x + b x + c x^3 + p b x^4 + p^2 a x^5 + p^3 x^6.
\]
Using {\em Magma}'s \verb@ZetaFunction@ again, and taking the first two good primes in each case, one gets
\begin{align*}
L_5(Q_0,x) & =  1 + x + 3x^2 + x^3 + \cdots, &
L_7(Q_0,x) & = 1 - x + 4x^2 - 11x^3 +\cdots, \\
L_5(Q_1,x) & = 1 + x + 3x^2 + x^3 + \cdots, &
L_7(Q_1,x) & =  1 - x + 8 x^2 - x^3 + \cdots, \\
L_5(Q_2,x) & = 1 + x + x^2 + 11x^3 + \cdots, &
L_7(Q_2,x) & = 1 - x + 8x^2 - x^3 + \cdots
\end{align*}
One has coincidences $L_5(Q_0,x) = L_5(Q_1,x)$ and
$L_7(Q_1,x) = L(Q_2,x)$, with the second polynomial
being reducible: $ (1-x+7x^2)(1+x^2 + 49 x^4)$.  
The generic behavior is that all three $L_p(Q_i,x)$
are different and their splitting fields
are disjoint extensions of $\Q$, each with
Galois group the wreath product $S_2 \wr S_3$ of order
$48$.  

The behavior of the curves here differs sharply from 
the behavior of the curves in \S\ref{UL}.   To describe 
this difference, we will use the language of motives, 
referring to the unconditional theory of \cite{And}.  
Note however, that the language of Jacobians
would suffice for the current comparison.  
Similarly, one could use the alternative language
of Artin representations for \S\ref{g2small}. 
But for uniformity, and certainly to include
the general case as represented by
\S\ref{g2big}, the language of motives is
best.  

The difference between they $Y_i$ of \S\ref{UL} and the
$\Q_i$ here goes as follows.  The two curves $Y_i$ from \S\ref{UL}
give rise to a single rank six motive
 $M = H^1(Y_2,\Q) \subset H^1(Y_1,\Q)$.  
Moreover the potential automorphism $(x,y) \mapsto (x,iy)$
causes the motivic Galois group of $M$ to be the ten-dimensional
conformal unitary group $CU_3.2$.  In contrast, the
motives $M_i = H^1(Q_i,\Q)$ here are all different,
as is clear from their different $L$-polynomials.  
Moreover, their motivic Galois groups 
are all as big as possible, the full $22$-dimensional
conformal symplectic group $CSp_6$.  

  While the different $L_p(Q_i,x) \in \Z[x]$ have very little to do with each other, 
their reductions to $\F_2[x]$ coincide, as illustrated with 
primes $5 \leq p \leq 97$:
\[
\begin{array}{l  lll}
\mbox{Class}(p) & \lambda_{28}(p) & L_p(Q_i,x) \in \F_2[x] &  \mbox{Primes } p \\
\hline
3B   &   3^9 1& (x+1)^2 \left(x^2+x+1\right)^2 & 89 \\
7AB   & 7^4&  \left(x^3+x+1\right) \left(x^3+x^2+1\right) & 5,13,29,53,61,73,97 \\
8AB   &  8^3 21^2 &  (x+1)^6 & 37,41 \\
12AB  & 12^2 3 1  &   \left(x^2+x+1\right)^3 &  17 \\
6b   &  6^4 3 1&    (x+1)^2 \left(x^2+x+1\right)^2 & 11,19 \\
8c  & 8^3 4&   (x+1)^6 & 43,67,79,83 \\
12c,12d & 12^2 3 1&    \left(x^2+x+1\right)^3 & 7,23,31,47,59,71. \\
   \end{array}
\]
This table shows very clearly how $S_0(-4,-3,z)$ functions
as a $2$-division polynomial.  All three $S_i$, 
arbitrarily specialized, similarly capture
the mod $2$ behavior of corresponding $L$-polynomials.

\section{$2$-division polynomials of Dettweiler-Reiter $G_2$ motives}
\label{G2}  
This section explains how the cover $\pi_1 : X_1 \rightarrow U_{3,1,1}$ is 
related to rigidity in the algebraic group $G_2$ in two 
ways.  The last subsection presents some sample
analytic calculations with $L$-functions.  

\subsection{Rigidity in general}
         In the mid 1990s, Katz \cite{Katz} developed a powerful theory of rigidity of tuples $(g_1,\dots,g_z)$ satisfying
$g_1 \cdots g_z = 1$ in ambient groups of the form $GL_n(E)$, with $E$ being an algebraically closed field.  
There is presently developing a theory of rigidity of tuples in $G(E)$ for other ambient algebraic groups $G$;
particularly relevant for us is \cite{DR}, where $G$ is either $G_2$ or $SO_7$.  
In general, if $G$ is simple modulo its finite center we say that a tuple $(C_1,\dots,C_z)$ is numerically rigid if 
\begin{equation}
\label{rigidnum}
\sum_{i=1}^z \ced_G(C_i) = (n-2) \dim(G).
\end{equation}
Here for $C_i$ a conjugacy class containing an element $g_i$, 
the integer $\ced_G(C_i) = \ced_G(g_i)$ is the dimension
of the centralizer of $g_i$ in $G(E)$.  

The Malle-Matzat case provides a convenient example in Katz's original context.  As explained in \cite[\S 8]{ABC},
after a quadratic base change the class triple $(4b,2b,12AB)$ becomes 
 $(12A,2A,12B)$ in $\Gamma = SU_3(\F_3)$.  Pushed forward to $SL_3(\overline{\F}_3)$, the classes $12A$ and $12B$ are regular and
 so have centralizer dimension $\mbox{rank}(SL_3) = 2$.  The class $2A$ is a reflection and has centralizer 
 $GL_2(\overline{\F}_3)$ with dimension $4$.  The rigidity condition \eqref{rigidnum} becomes $2+4+2=1 \cdot  8$ and is thus satisfied.  

\subsection{Groups $G_2(2)$ and $G_2$-rigidity} 
\begin{figure}[htb]
\includegraphics[width=4.8in]{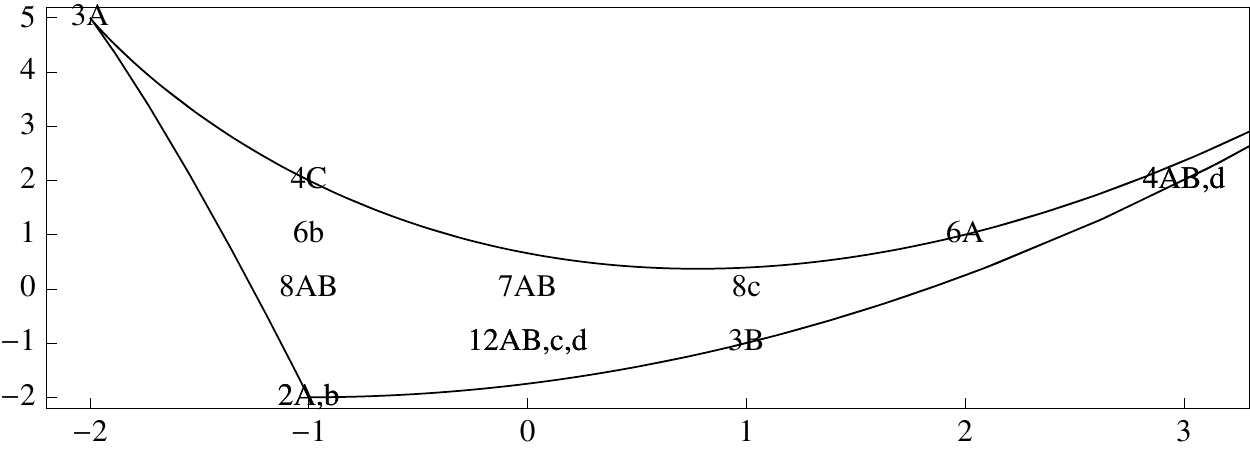}
\caption{\label{g2classpict} The image of the class-set $G_2(2)^\natural$ inside the class space $G_2^{c \natural}$.}
\end{figure}

Our group $\Gamma.2 = G_2(2)$ embeds into the fourteen-dimensional compact Lie
group $G^c_2$.  Figure~\ref{g2classpict}
illustrates the associated map $G_2(2)^\natural \rightarrow G_2^{c \natural}$
 on the level of conjugacy classes, which
is no longer injective.  The fundamental characters $\chi$ and $\phi$ 
of $G_2$ have degrees $7$ and $14$ respectively, and
the set $G_2^{c \natural}$ becomes the indicated triangular region 
in the $\chi$-$\phi$ plane.    The unique class in $G_2(2)^\natural$
which is outside the window is the identity class $1A$ 
at the point $(\chi,\phi)=(7,14)$.  
 
 We have drawn Figure~\ref{g2classpict}  to facilitate the analysis of rigidity in $G_2(\C)$.     First consider
 classes which intersect the compact group $G_2^c$, and are thus represented by points in 
 the closed triangular region drawn in the figure.  If $g$ represents one of the vertex classes labeled by $1A$, $3A$, and $2A,b$ respectively, 
 its centralizer has type $G_2$, $SL_3$, and $SL_2 \times SL_2$, thus dimension $14$, $8$, and $6$ 
 respectively.  For $g$ representing a class otherwise on the boundary, the centralizer has type $GL_2$
 and hence dimension $4$.  For $g$ in the interior, the class is regular and so the centralizer 
 dimension is $\mbox{rank}(G_2)=2$.  For general semisimple elements in $G_2(\C)$ the situation is the same:
 centralizer dimensions are $14$, $8$, and $6$ for the three special classes
 already considered, $4$ for classes on the algebraic curve 
  corresponding to the boundary, 
 and $2$ otherwise.
 
 \subsection{$G_2$-rigidity of $(3A,3A,3A,4B)$} 
 \label{g2small} The first-listed quadruple for $\pi_1$ in Prop.~\ref{rigidprop} is 
 $(3A,3A,3A,4B)$.    Using the determinations associated to Figure~\ref{g2classpict},
 the left side of \eqref{rigidnum} becomes $8+8+8+4=28$ which agrees with the right side $(4-2) 14 = 28$.  
  Thus $(3A,3A,3A,4B)$ is $G_2$-rigid.  We are 
 not pursuing this connection here, but it seems possible to write down a corresponding rank seven
 differential equation with finite monodromy.   From the algebraic solutions to this differential
 equation, one could perhaps construct the cover $X_1 \rightarrow U_{3,1,1}$ in a third way.
 
 \subsection{Orthogonal rigidity of a lift of $(2A,2A,3A,4A)$}  
 \label{g2big} The last-listed quadruple for $\pi_1$ in 
 Prop.~\ref{rigidprop} is $(2A,2A,3A,4A)$.   This tuple fails to be $G_2$-rigid, as now
 the left side of \eqref{rigidnum} is $6+6+8+3 = 24$ which is less than $28$.    However one does have
 rigidity of a lift as follows.  

   The following matrices were sent to me by Stefan Reiter in July 2013.  
   \begin{align*}
a & =   \left(
                 \begin{array}{rrrrrrr}
                  1 &  & & &&&  \\
                  & 1 &&  & & &  \\
                   &  & 1 & & & & \\
                  -3 & 1 & & 1 & & & \\
                  3 & -1 &  & & 1 & &  \\
                  9 & -3 & & &  & 1 & \\
                  -1 &  & 3 & -1 & 2 & -1 & 1 \\
                 \end{array}
                 \right)  
                 & \sim 
                 \left(
                 \begin{array}{rrrrrrr} 
                 1&1&&&&&\\
                 &1&&&&&\\
                 &&1&1&&&\\
                 &&&1&&&\\
                 &&&&1&&\\
                 &&&&&1&\\
                 &&&&&&1\\
                 \end{array}
                 \right) \\
                 &\\
b & =  \left(
                 \begin{array}{rrrrrrr}
                  1 &  &  & & 3 & -1 &  \\
                  & 1 &  &  & 9 & -3 & \\
                  &  & -2 & 1 &  &  &  \\
                  & & -9 & 4 & & & \\
                   &  &  &  & 1 & & \\
                 & & &  &  & 1 &  \\
                   &  & -3 & 1 &&  & 1 \\
                 \end{array}
                 \right) 
                  & \sim 
                 \left(
                 \begin{array}{rrrrrrr} 
                 1&1&&&&&\\
                 &1&&&&&\\
                 &&1&1&&&\\
                 &&&1&&&\\
                 &&&&1&&\\
                 &&&&&1&\\
                 &&&&&&1\\
                 \end{array}
                 \right) \\
                 &\\
 c & = 
                 \left(
                   \begin{array}{rrrrrrr}
                   1 & -1 & & & &  & -3 \\
                   3 & -2 &  & &  &  &  \\
                    & & 1 & -1 &  &  & 3 \\
                    & & 3 & -2 &  & & 6 \\
                    &  & &  & 1 & -1 & -3 \\
                    &  & &  & 3 & -2 & \\
                   &  & & &  &  & 1 \\
                  \end{array}  
                  \right) 
                                  & \sim 
                 \left(
                 \begin{array}{rrrrrrr} 
                 \! \omega  &&&&&&\\
                 &\! \omega&&&&&\\
                 &&\! \omega&&&&\\
                 &&& \! \overline{\omega}&&&\\
                 &&&&\! \overline{\omega}&&\\
                 &&&&& \! \overline{\omega}&\\
                 &&&&&&1\\
                 \end{array}
                 \right) \\
                 &\\
d & =                   \left(
                  \begin{array}{rrrrrrr}
                   10 & -5 & & & 9 & -5 & -6 \\
                   15 & -8 & &  & 18 & -9 & -9 \\
                    &  & 1 & &  &  &  \\
                   -3 & 4 & -3 & 1 & -6 & 3 & 3 \\
                   9 & -5 & & & 10 & -5 & -6 \\
                   18 & -9 &&& 15 & -8 & -9 \\
                   -2 & 1 && & -2 & 1 & 1 \\
                  \end{array}
                  \right)
                   & \sim 
                 \left(
                 \begin{array}{rrrrrrr} 
                 1&1&&&&&\\
                 &1&1&&&&\\
                 &&1&&&&\\
                 &&&1&1&&\\
                 &&&&1&1&\\
                 &&&&&1&\\
                 &&&&&&1\\
                 \end{array}
                 \right) 
\end{align*}
These matrices satisfy $abcd=1$ and they generate a subgroup
of $GL_7(\C)$ with Zariski closure of the form $G_2(\C)$.  
On the one hand, reduced to $GL_7(\F_2)$, these 
matrices generate a copy of $G_2(2)'$ with
$a$, $b$, $c$, and $d$ respectively in $2A$, $2A$, $3A$ and $4A$.  
On the the other hand, considered in $GL_7(\C)$, the matrices have Jordan
canonical forms as listed on the right, with $\omega = \exp(2 \pi i/3)$. 

Consider $a,b,c,d \in G_2(\C) \subset SO_7(\C) \subset SL_7(\C)$.  Centralizer dimensions are calculated in \cite[\S3]{DR} 
 and the numerics
associated with \eqref{rigidnum} are as follows.
{\renewcommand{\arraycolsep}{1pt} 
\[
\begin{array}{c|ccccccccccccc}
G & \ced_G(a) &+& \ced_G(b) &+& \ced_G(c) &+& \ced_G(d) &\;\;\; & \;\;\;\;\;\;\;\;\;\;\;& \;\;\; &    2 \dim(G) \\
\hline
G_2 & 8 & + & 8 & + & 8 & + & 4& = & 28 &  = & 28 \\ 
SO_7 & 13 & + & 13 & + & 9 & + & 7& = & 42 &  = & 42 \\ 
SL_7 & 28 & + & 28 & + & 28 & + & 16 & = & 90 & < & 96  
\end{array}
\]
}

\noindent Thus the quadruple $([a],[b],[c],[d])$ is $G_2(\C)$- and $SO_7(\C)$-rigid.  
However it is not $SL_7(\C)$-rigid, and so does not fit into Katz's 
original framework. 

Dettweiler and Reiter classify tuples of classes in $G_2(\C)$ which are $SO_7(\C)$ rigid
in \cite{DR}.   Thus $([a],[b],[c],[d])$ is in their classification.  In fact, it appears as the first line of the 
table in \S5.4.     Being $SO_7(\C)$-rigid is a stronger condition than being $G_2(\C)$-rigid. 
It implies from \cite{DR} that there is a corresponding rank seven motive over $\Q(p,q)$ with motivic Galois group $G_2$. 

\subsection{Division polynomials and $L$-functions} 
\label{GL}  In \S\ref{UL} and \S\ref{SL} we have discussed $L$-polynomials $L_p(M,x)$ for 
certain motives $M = H^1(\mbox{curve},\Q)$.   Putting these $L$-polynomials together,
including also $L$-polynomials at bad primes, one gets a global $L$-function
\begin{equation}
\label{GlobalL} L(M,s) = \prod_{p} L_p(M,p^{-s})^{-1}.
\end{equation}
This $L$-function is expected to have standard analytic properties, including an analytic
continuation and a functional equation with respect to $s \leftrightarrow 2-s$.  Normalizing the motives from \S\ref{g2small} and \S\ref{g2big} 
to have weight $0$, one likewise expects good analytic properties of corresponding $L(M,s)$, involving now functional 
equations $s \leftrightarrow 1-s$. 

We do not know yet how to compute $L$-polynomials in the context of \S\ref{g2big}, where
the motivic Galois group is generically the fourteen-dimensional algebraic group $G_2$.   
However the computation of $L$-polynomials is feasible in the setting of \S\ref{g2small} where 
the motivic Galois group is just the finite group $G_2(2)$.   In fact, as commented already in
\S\ref{SL}, we are using motivic language mainly because it is the
natural general context for division polynomials.  The particular motives from \S\ref{g2small} correspond
to finite-image Galois representations and so this language could 
be avoided.  

In the setting of Section~\ref{DM}, Section~\ref{shioda}, and \S\ref{g2small}, analytic computations
with global $L$-functions \eqref{GlobalL} are possible on a numerical level.   To illustrate this, 
we consider the motive $M$ from \S\ref{g2small} associated to the specialization
point used in \S\ref{UL} and \S\ref{SL}, namely
$(u,v) = (-4,-3)$.   This motive corresponds to 
the seven-dimensional irreducible representation of $G_2(2)$ into $SO(7)$.  It is natural
here to twist by the Dirichlet character $\chi$ given on odd primes $p$ by
$\chi(p) = (-1)^{(p-1)/2}$.   The twisted motive
 $M'$ corresponding to the other seven-dimensional irreducible
representation of $G_2(2)$.   At the level of good $L$-polynomials, passing back
and forth between $M$ and $M'$ means replacing $x$ by $\chi(p) x$.

Let $p \geq 5$ be a prime.  The corresponding Frobenius class $\Fr_p$ can 
usually be deduced from Table~\ref{wtable} 
from the mod $p$ factorization partition of $S_0(-4,-3,z)$ and 
the class of $p$ modulo $4$.  To make the necessary distinction between
$3A$ and $3B$, we use the factorization partition of the resolvent $f_{36}(4,x)$ presented in \eqref{f36}.  
The $(\chi,\phi)$-coordinates of $\Fr_p$ on Figure~\ref{g2classpict} then yield the $L$-polynomial
\[
L_p(M,x) = 1 - a x + b x^2 - c x^3 + c x^4 - b x^5 + a x^6 - x^7.
\]
Here $a = \chi$,  $b = \chi+\phi$, and $c = a + a^2-b$.    

The necessary $2$-adic and $3$-adic analysis 
for obtaining conductors and bad $L$-polynomials 
is begun in Prop.~\ref{badanalysis} below.  For $L(M,s)$
the conductor is $2^{20} 3^{12}$,   the decomposition
of the exponents as a sum of seven slopes being as follows.
\begin{align*}
\mbox{At 2: } \; 20 & =  6 \cdot 3 +2. &
\mbox{At 3: } \; 12 & =  6 \cdot \frac{11}{6} + 1. 
\end{align*}
Since all slopes are positive, the bad $L$-polynomials are 
$L_2(M,x) = L_3(M,x)=1$.   For $L(M',s)$,
the slopes are all the same except the $2$-adic slope
$2$ is now $0$, so that the conductor drops to $2^{18} 3^{12}$. 
 Slopes of $0$ contribute to the degree of $L$-polynomials, and
in this case $L_2(M',x)=1-x$ while still 
$L_3(M',x)=1$.

\begin{figure}[htb]
\includegraphics[width=4.5in]{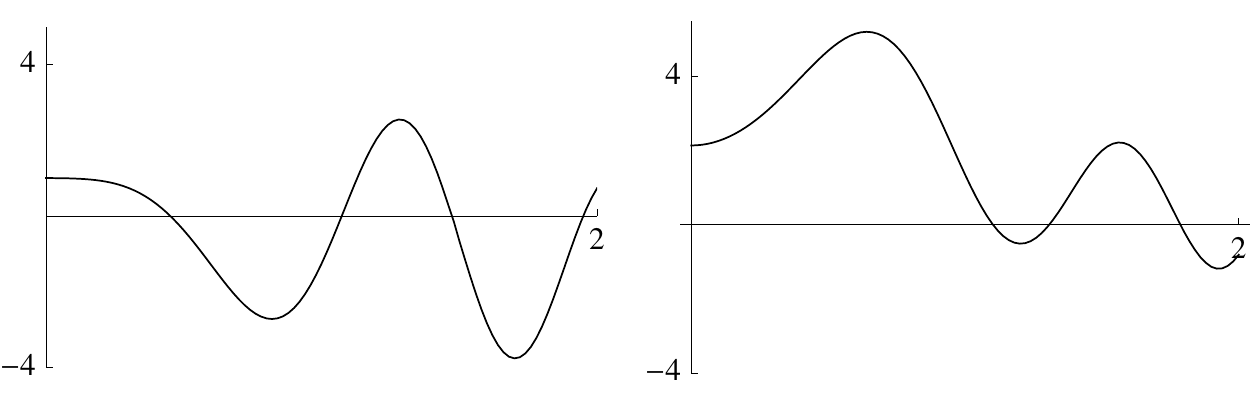}
\caption{\label{critline} Graphs of $L^*(M,\frac{1}{2}+i t)$ (left) and 
$L^*(M',\frac{1}{2}+i t)$ (right) } 
\end{figure}

In principle, {\em Magma}'s Artin representation and $L$-function packages \cite{Mag}, 
both due to Tim Dokchister, should 
do all the above automatically, given
simply $S_0(-4,-3,z)$ as input.   However the inertia groups
at $2$ and $3$ are currently too large, and so {\em Magma} 
can only be used with the above extra information 
at the bad primes.  It then outputs numerical values for arbitrary 
$s$, on the assumption that standard conjectures hold.  
Particularly interesting $s$ include those of the form 
$\frac{1}{2} + i t $ with $t$ real, i.e.\ those on the critical
line.  Here one multiplies $L$ by a phase factor depending
analytically on $t$ to obtain a new function $L^*$ taking
real values only.   Figure~\ref{critline} presents plots for
our two cases, numerically identifying zeros on the
critical line.

To obtain analogous plots of $L^*(M,\frac{w+1}{2}+it)$ for
a general weight $w$ motive, such as the 
weight one motives from  \S\ref{UL} and 
\S\ref{SL}, division polynomials do not at all
suffice.    Here one needs the much more complete
information obtained from point counts,
like the $L_p(M,x)$ presented in \S\ref{UL} and
\S\ref{SL} for $p=5$ and $p=7$.   However
division polynomials can still be of assistance
in obtaining the needed information 
at the bad primes.

\section{Specialization to three-point covers}
\label{threepoint}
       In \S\ref{basecurves} we find projective lines $P$ in $\overline{U}_{3,1,1}$ and
 $\overline{U}_{3,2}$ suitably intersecting the discriminant 
 locus in only three points.  In \S\ref{3point} we consider the covers 
 obtained by the preimages under $\pi_1$ and $\pi_2$ of these lines.  We 
 thereby construct some of the three-point covers $X_P \rightarrow P$
 mentioned in \S\ref{fivecovers}.  As stated previously, 
  it would be hard to construct these covers directly because these $X_P$ 
 always have positive genus.   In \S\ref{recovering} we apply 
 quadratic descent twice to a cover 
 $X_P \rightarrow P$ coming from a curve 
 $P \subset U_{3,1,1}$ and recover the Malle-Matzat cover
 \eqref{mm}.   
 
\subsection{Curves in $\overline{U}_{3,1,1}$ and $\overline{U}_{3,2}$} 
\label{basecurves} 
\begin{figure}[t]
\includegraphics[width=4.5in]{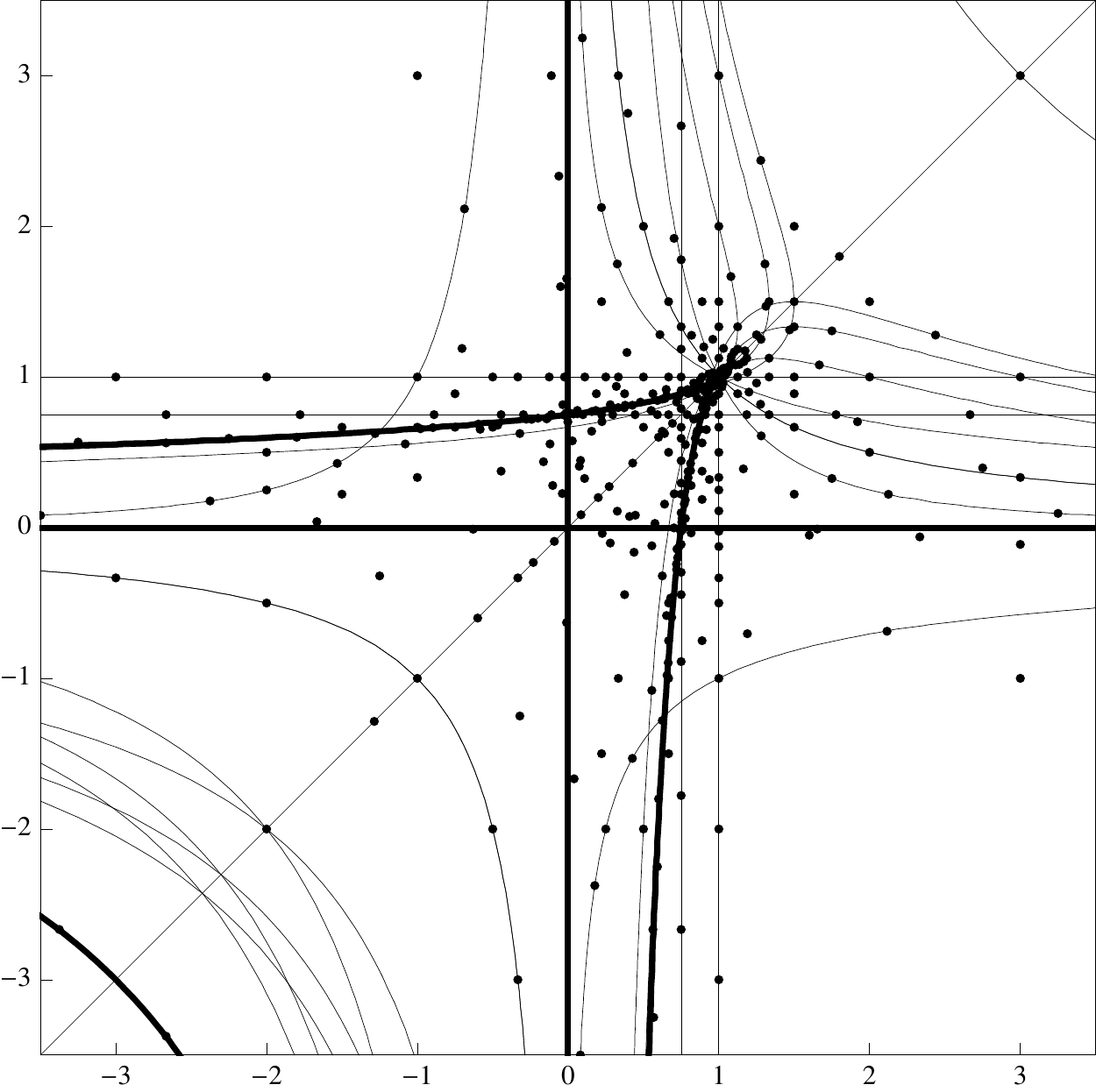}
\[
\begin{array}{lc|cc}
\multicolumn{2}{l|}{\;} &  p & q  \\
     \hline
G'& \bullet &  t & 1 \\
H'& \bullet &   {8 t}/{9} & {3}/{4} \\
I'&&   (9 t-5)/4 & {1}/{p^2}  \\
J'& &   {(t+3)}/{4} & {(4 p-3)/p^2}  \\
K'&a&   (t+2)/3 & (3 p-2)/p^2\\
L'&e &   3 (t+4)/16 & \!\! 9(4 p-3)/8 p^2 \!\! \\
    F^* & \bullet f &     {3}/{(4 t-1)} & {3}/(4 t-1)  \\
M^*&&   t & 1/t \\
 B^*&b&    -3 t & -{3}/{t}  
   \end{array}
   \;\;\;\;\;\;
\begin{array}{cc|cc}
\multicolumn{2}{c|}{\;} &  p & q  \\
    \hline
G''& \bullet  &  1 & t\\
H''& \bullet &    3/4& 8t/9 \\
I''&& {1/q^2} & (9 t-5) /4\\
J''&&    {(4 q-3)}/{q^2} & {(t+3)}/{4}  \\
K''& c&  (3 q-2)/{q^2} & (t+2)/3  \\
L''& d&  \!\! 9 (4 q-3)/8 q^2 \!\! & 3 (t+4)/16 \\
\multicolumn{4}{c}{\;} \\
\multicolumn{4}{c}{\;} \\
\multicolumn{4}{c}{\;} \\
   \end{array}
   \]
\caption{\label{u311pict} Top: The $p$-$q$ plane $U_{3,1,1}(\R)$.  Discriminant loci (thick), 
bases of three-point covers (thin), and specialization points are drawn in.  Bottom: 
parametrizations of the bases for three-point covers}
\end{figure}

The top half of 
Figure~\ref{u311pict} is a window on
the real points of the naive completion
$\overline{U}_{3,1,1}= \bbP^1_p \times \bbP^1_q$. 
The discriminant locus ${Z}_{3,1,1}$ consists of the two coordinate
axes, the two lines at infinity, and the solution curve $D_1$ of
\[
p^2 q^2 - 6pq + 4 p + 4 q - 3=0.
\]
The five lightly 
drawn straight lines intersect ${Z}_{3,1,1}$ in just three points, not counting
multiplicities.       The ten other lightly drawn curves have
the same three-point property, although it is not visually evident.  
The points drawn in Figure~\ref{u311pict} will be discussed in the next section.

The bottom half of Figure~\ref{u311pict} names and parametrizes the fifteen lightly drawn
curves in the top half.   Each name is a superscripted letter.  The five bulleted 
curves are the straight lines.  There are other natural coordinate systems on the $p$-$q$-plane,
and each of the other curves appears as a line in at least one
of these coordinate systems.  We are emphasizing the coordinates
$p$ and $q$ because they make the natural involution of $U_{3,1,1}$
completely evident as $p \leftrightarrow q$.   The three 
curves labeled $T^*$ are stable under this involution.
The remaining twelve curves form six interchanged pairs:  
 $T' \leftrightarrow T''$.   
Six of the fifteen curves are images of lines in the cubic cover $U$.  These
source lines in $U$ are indicated by $a$, $b$, $c$, $d$, $e$, and 
$f$.   

\begin{figure}[t]
\includegraphics[width=4.5in]{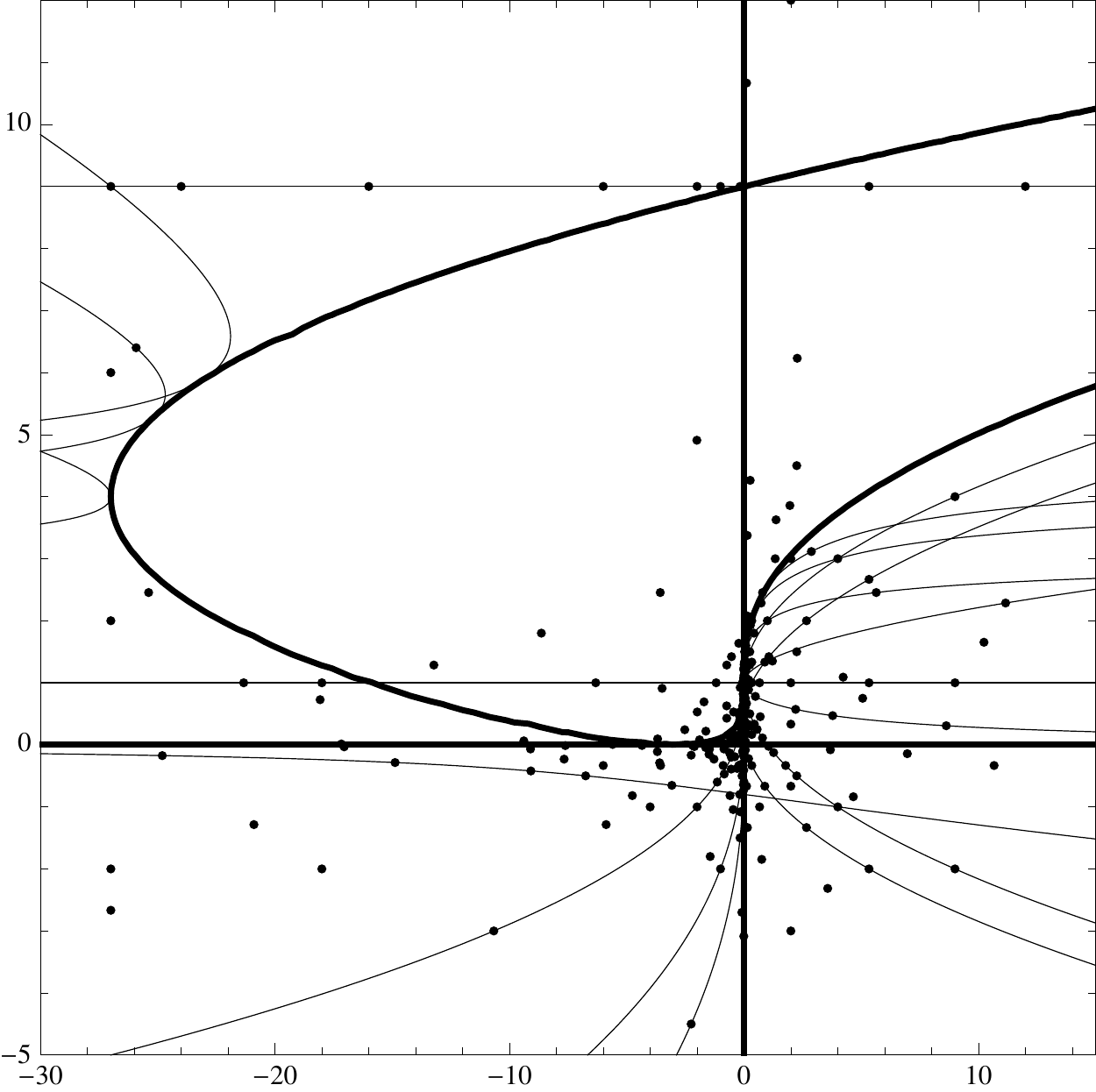}
{\renewcommand{\arraycolsep}{1.3pt}
{
\[
 \begin{array}{l|cc}
& a & b  \\
 \hline
 A & -{27 (b-1)^3}/{(5 b-9)} & 9 (t-1)/{(5 t-9)}  \\
B\bullet & -48 (t-1) & 9  \\
C  & - (b-9) (5 b-9)/27 & 9 (t-1)/{(5 t-1)} \\
E    & -{(b-1)^3}/{(b-3)} & {3 (t-1)}/{(t-3)}  \\
F  & -{3}/{\left(9 t^2-7 t+1\right)^2} & \frac{9 (t-1) t}{9 t^2-7 t+1}  \\
M\bullet &     -16 t  &  1         
\end{array}
\;\;\;\;\;\;
\begin{array}{l|cc}
 &  a & b  \\
 \hline
G  &       (b-1)^2 & t  \\
H  &       b (3 b-2)/3 & {2 t}/{3}  \\
I  &    {(b-1)^2 (5 b+4)/b} & {4}/{(9 t-5)}  \\
J &    -{(b-1)^2 b/3 (b-4)} & {4 t}/{(t+3)}  \\
K  &    -{(b-1)^3}/{(b-3)} & {3 t}/{(t+2)}  \\
L &     \frac{-b^2 (10 b-9)}{27 (2 b-9)} & {9 t}/{2 (t+4)}  \\
\end{array}
\]
}
}
\caption{\label{pictab} Top: The $a$-$b$ plane $U_{3,2}(\R)$.  Discriminant loci (thick), 
bases of three-point covers (thin), and specialization points are drawn in.  Bottom: 
parametrizations of the bases for three-point covers.}
\end{figure}

Figure~\ref{pictab} is the analog of Figure~\ref{u311pict} for $\overline{U}_{3,2} = \bbP^2_{a,b}$
and we will describe it more briefly, focusing on differences. 
The discriminant locus $Z_{3,2}$  has four components, the two coordinate
axes, the line at infinity, and the curve $D_2$ with equation
\[
a^2-2 a b^2+12 a b+6 a+b^4-12 b^3+30 b^2-28 b+9=0.
\]
The light curves each intersect the discriminant locus in three points, where
this time a contact point with $D_2$ does not count if the local intersection
number is even.  Despite the relaxing of the three-point condition, we have found
only twelve such curves.  The five curves $A$, $B$, $C$, $E$, and $F$ 
are images of generically bijective maps from 
curves $a$, $b$, $c$, $e$, and $f$ in $U$.   Curve $d$ in 
$U$ double covers $B$, and so does not 
have its own entry on Figure~\ref{pictab}.    For $T = G$, $H$, $I$, $J$, $K$, and $L$,
the curve $T \subset U_{3,2}$ comes from $T'$ and $T''$ in $U_{3,1,1}$ via \eqref{ab}.  
Finally $M \subset U_{3,2}$ is double-covered by $M^*$ in $U_{3,1,1}$.

\subsection{Three-point covers with Galois group $\Gamma.2$}  
\label{3point}
The previous subsection concerned the base varieties $U_{3,1,1}$ 
and $U_{3,2}$ only.  For quite general covers $X \rightarrow U_\nu$,  one
gets three-point covers $X_P \rightarrow P$ by specialization
to the $P \subset X_\nu$ listed there.   We now apply this
theory to our particular covers $\pi_1 : X_1 \rightarrow U_{3,1,1}$
and $\pi_2 : X_2 \rightarrow U_{3,2}$.   Because of the explicit
parametrizations in Figures~\ref{u311pict} and \ref{pictab}, our
bases are now coordinatized projective lines $\bbP^1 = \bbP^1_t$.  

\begin{table}[htb]
\[
\begin{array}{ccc|ccc|cr|cc}
X_0 & X_{311} & X_{32}  & C_0 & C_1 & C_\infty & g_{28} & g_{36} & \bar{\mu} & \mu \\
\hline    
        &      H''    &      &  4A    &   4B &   3B  &  -   &  - & 0.\overline{3}  &  0   \\
      &    I''        &      &   4A &  12A   &   2A  &  -   &    - &  0.\overline{3}  &  0    \\
\hline
      b     &      B^*          & B & 6A & 2A & 8A & 1 & 0 &1 & 1  \\
 \hline
       &              &   M & 12A & 2A & 8B & 2 & 2 & 1 & 1 \\
       &              & G  &    4A  &  6A    &  3B   &   2  &  2  &   1 & 1    \\
      &    H', G''   &      & 12A & 4A  &3B & 2 & 5 & 1  &  1    \\
 \hline
      e     &     L'         & E,K & 4C & 4A & 8A & 3 & 3 & 1 & 1  \\
  &       G'      &  H & 3A & 12A & 3B & 3 & 5 &1 & 1 \\
\hline
a     &      K'        & A & 4A & 8A & 8B & 4 & 7 & 1   & 1 \\
c     &      K''         & C,I  & 3A & 8A & 6A & 4 & 6 &  1 & 1 \\
d     &     L''         &  & 6A & 4A & 6A & 4 & 5 &  1  & 1  \\
\hline
f     &      F^*,I'     & F  & 4A & 8B  & 12B  & 5 & 8 & 1 & 1 \\ 
      &       J'     &      &  4A  &  12A  &  8B &  5  &  8  &  1  &   1  \\
            &              &   L &  12A    &   3A  &  8A   &   5  &   8  &  1  &  1   \\
\hline
                 &       M^*    & J   & 6A & 12A & 8B & 7 & 10  & 5 & 5   \\
\hline
      &      J''     &      &  12A   &   12A   &   6A  &   8  &  11   &  4.08\overline{3}  &  3   \\
    \end{array}
\]
\caption{\label{sixteen} Sixteen three-point covers obtained from $\pi_1$ and $\pi_2$ 
by specialization }
\end{table}

Table~\ref{sixteen} gives the results.  The first two lines illustrate
the general phenomenon where Galois groups sometimes 
become smaller under specialization.  Here the covers have  
Galois groups of order $216$ and $432$ respectively, thus
of index $56$ and $28$ in $\Gamma.2$.    The covers
$X_{28} \rightarrow \bbP^1$ each split into a 
genus one cover $X_{27} \rightarrow \bbP^1$ and
the trivial cover $\bbP^1 \rightarrow \bbP^1$.  

The next fourteen lines each give a cover
$X_{28} \rightarrow \bbP^1$ with Galois group 
all of $\Gamma.2$.   They are sorted by the 
genus $g_{28}$ of this cover.   In
most cases, more than one base curve 
$\bbP^1$ yield isomorphic covers, after
suitable permutations of the three cusps $\{0,1,\infty\}$.  
The local monodromy classes in $\Gamma$ 
always correspond to the first-listed parametrized base. 
These classes are unambiguously
determined, except for a simultaneous interchange
$4A \leftrightarrow 4B$, $8A \leftrightarrow 8B$,
$12A \leftrightarrow 12B$, coming from the
outer automorphism of $\Gamma$.  We 
always normalize by making the first-listed 
interchanged class have an $A$ is its name.

Thus for example, specializing $S_1(p,q,x)$ at $(p,q) = (-3t,-3/t)$ from the
$B^*$ line of Figure~\ref{u311pict}, one gets
a polynomial in $\Z[t,x]$ with $554$ terms.   The local 
monodromy partitions
are $(6A,2A,8A)$ as printed.  Alternatively, specializing
$S_2(a,b,x)$ at $(a,b) = (-48(t-1),9)$ from
the $B$ line of Figure~\ref{pictab}, one gets a polynomal
in $\Z[t,x]$ now with $252$ terms.  The
monodromy partitions are the same, except
for the reordering $(C_0,C_1,C_\infty) = (2A,6A,8A)$.

Having specialized from two parameters down to one, it is now much more reasonable to print
polynomials giving equations $f_{28}(t,x)=0$ and $f_{36}(t,x)=0$ corresponding to the covers
in any of the last fourteen lines of Table~\ref{sixteen}.   We do this only in the case where genera are the smallest,
namely the third line:  
\begin{eqnarray}
\nonumber \lefteqn{
f_{28}(t,x) =} \\
\nonumber&&-t  \left(3 x^4-252 x^3+222 x^2-692 x-5\right) \cdot \\ 
\nonumber&&  \qquad  \left(81 x^{12}+2106 x^{11}+26001
    x^{10}+73332 x^9+268515 x^8+574938 x^7 \right. \\
\label{f28}   && \qquad \left. +618759 x^6+400896 x^5+184140 x^4+52752 x^3+8952 x^2+576
    x-32\right)^2 \\
 \nonumber   &&+2^{10} (1-t)(4 x+1) \left(9 x^4+18 x^3+48 x^2+18 x+1\right)^6 \\
 \nonumber   &&+ 3^9 (1-t) t (x-2)^8 x^2
    \left(x^2+8\right) \left(x^2-2 x-1\right)^8, \\
 \nonumber   &&\\
\nonumber \lefteqn{
f_{36}(t,x) =} \\
\nonumber && \qquad  \left(4 x^4-3\right)^3 \left(4 x^4-12 x^2+12 x-3\right)^6  \\ 
\label{f36} && \qquad -3^9 t (x-1)^4 \left(2
    x^2-1\right)^8 \left(2 x^2-2 x+1\right)^4.
\end{eqnarray}
Here the genera, namely $(g_{28},g_{36}) = (1,0)$ are the 
reverse of those of the Malle-Matzat covers.

 \subsection{Recovering the Malle-Matzat curve}  
 \label{recovering}
 The Malle-Matzat cover can be constructed from the last line of Table~\ref{sixteen} via
 two quadratic descents as follows.  The given cover $X_1 \rightarrow \bbP_t^1$ has
 ramification invariants $(12A,12A,6A)$.  Quotienting out by
 the involution $t \leftrightarrow 1-t$ on the base and its
 unique lift to $X_1$, one gets the descended cover $X_2 \rightarrow \bbP_s^1$,
 with $s = 4 t(1-t)$.    The ramification invariants of this cover are $(12A,2A,12B)$.
 Quotienting now by $s \leftrightarrow 1/s$ on the base and its unique 
 lift to $X_2$, one gets the twice descended cover $X_3 \rightarrow \bbP^1_u$,
 with $u = -(s-1)^2/4s$.   The ramification invariants of this cover are
 $(4b,2b,12AB)$, showing that it is the Malle-Matzat cover.  
 
 In other words, 
 \[
 m \left( \frac{(2 t-1)^4}{16 (t-1) t}  ,x \right)  \mbox{ and } S_1\left( \frac{16 t}{(t+3)^2},\frac{t+3}{4},z \right)
 \]
are two different polynomials defining the same degree 28 extension of $\Q(t)$.  The
left one is a quartic base-change of the Malle-Matzat polynomial $m(u,x)$ 
while the right is a specialization of $S_1(p,q,x)$.

\section{Specialization to number fields}
\label{numbfield} 
   In this final section, we discuss specialization to number fields
with discriminant of the form $2^j 3^k$.  
\S\ref{specpi} discusses fields obtained by specializing the 
$\pi_i$. \S\ref{specsumm} continues this discussion, 
involving also similar fields from other sources.  
\S\ref{analysis} discusses analysis of ramification
in general, with a 
field having Galois group $PGL_2(7)$ serving as 
an example.  \S\ref{spechard} concludes by 
analyzing the ramification of a particularly 
interesting field with Galois group $SU_3(3).2 \cong G_2(2)$.

\subsection{Specializing the covers $\pi_i$}
\label{specpi}
    In this subsection, we restrict attention to number fields with
 Galois group $\Gamma.2$ and discriminant of the form $2^j 3^k$.  
 Consider first the cover $X_0 \rightarrow U$.  We have found 216 ordered pairs $(u,v)$
 such that the corresponding number field $\Q[x]/F_0(u,v,x)$ has Galois
 group $\Gamma.2$ and discriminant of the form $2^j 3^k$.   Different
 specialization points can give isomorphic fields, and we found 
 147 number fields in this process.
 
 Next consider the covers $\pi_1 : X_1 \rightarrow U_{3,1,1}$ and
 $\pi_2 : X_2 \rightarrow U_{3,2}$.  Beyond images of 
 specialization points in $U(\Q)$, we found  248 pairs $(p,q)$ 
 and 177 pairs $(a,b)$ giving fields with  Galois
  group $\Gamma.2$ and discriminant of the form $2^j 3^k$.  
 We obtained 62 new fields arising from both covers,
 95 new fields arising from $\pi_1$ only, and 72 new fields
 arising from $\pi_2$ only.   Thus we found in total
 376 fields with Galois group $\Gamma.2$ and discriminant of the form
 $2^j 3^k$.  
 
 \begin{figure}[htb]
\includegraphics[width=4.8in]{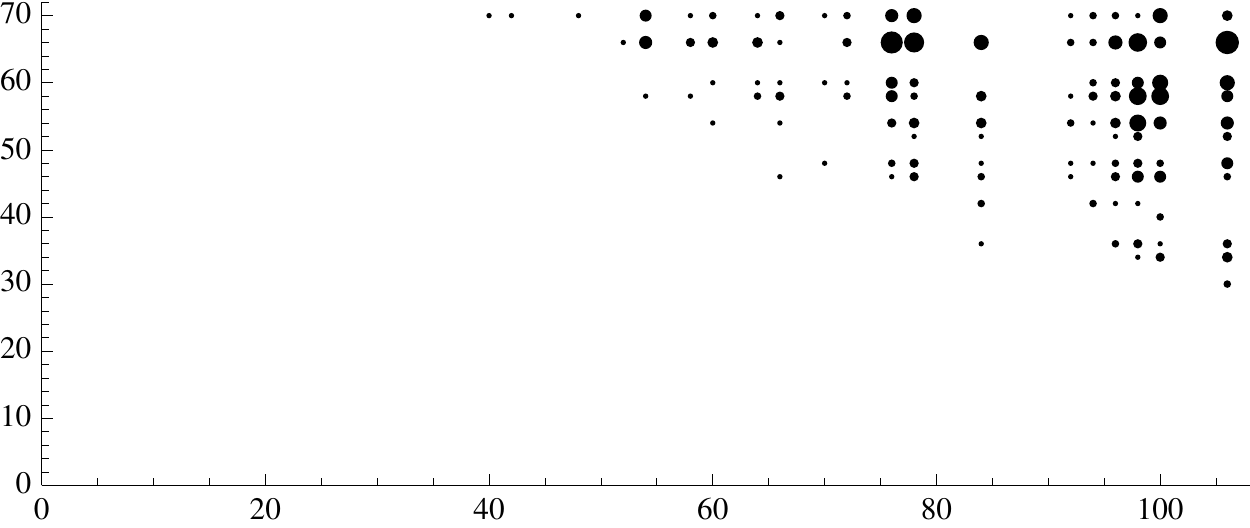}
\caption{\label{g2discfs} Pairs $(j,k)$ arising from field discriminants $2^j 3^k$ from 
specializations of $F_1(p,q,x)$ and $F_2(a,b,x)$}
\end{figure}

Figure~\ref{g2discfs} indicates the pairs $(j,k)$ arising from field discriminants $2^j 3^k$ 
of one of these 376 fields.  The area of the disk at $(j,k)$ is proportional to the
number of fields giving rise to $(j,k)$.  In $36$ cases, this field is unique.
The largest multiplicity is $19$, arising from $(j,k) = (106,66)$.  
The smallest discriminant is $2^{66} 3^{46}$, coming from just
one field.  This field arises from eight sources,
\begin{align}
\nonumber (u,v)& = (-4,-3), (-\frac{1}{2},1), (\frac{1}{2},3), (4,-3),(-32,1), (-\frac{32}{81},\frac{49}{81}),  \\
\label{repetition} (p,q)& = (1, \frac{1}{2}), \\
\nonumber (a,b)& = (-\frac{27}{4},-\frac{1}{2}).
\end{align}
The largest discriminant $2^{106} 3^{70}$ arises from four fields.  

The phenomenon of several specialization points giving rise to a single field
is quite common in our collection of covers $\pi_i$.   The octet in \eqref{repetition}
is the most extreme instance, but there are many other multiplets
as altogether
$216+248+177=641$ different specialization points give rise to
only 376 fields.  This repetition phenomenon 
 is discussed for a different cover in \cite[\S6]{ABC}, where it is explained by a 
 Hecke operator.   It would be of interest to give 
 a similar automorphic explanation of the very large
  drop $641 \rightarrow 376$.   Ideally, such a description
  would follow through on one of the main points of 
  view of Deligne and Mostow \cite{DM1,DM}, by describing
  all our surfaces via uniformization by the unit ball in 
  $\C^2$.

\subsection{Summary of known fields}
\label{specsumm}
We specialized the Malle-Matzat cover in \cite[\S8]{ABC} to obtain fields
with discriminant of the form $2^j 3^k$.    From 
$t = 1/2$ we obtained a field with Galois group $\Gamma$, while 
from $41$ other $t$ we obtained $41$ other fields with Galois group 
$\Gamma.2$.    While in our covers $\pi_i$ the $.2$ always corresponds
to the quadratic field $\Q(i)$, in the Malle-Matzat cover 
general $\Q(\sqrt{\partial})$ arise.  

Sorting all the known fields by $\partial \in \Q^\times/\Q^{\times 2}$, including two additional fields from \cite{CoversM05} with
$\partial=2$ and $\partial=6$,  one has the following result.
\begin{Proposition}
\label{g2fields}  There are at least $409$ degree twenty-eight fields
with Galois group $\Gamma$ or $\Gamma.2$ and discriminant of the form
$2^j 3^k$.  Sorted by the associated quadratic algebra $\Q[x]/(x^2-\partial)$, 
these lower bounds are 
\[
 \begin{array}{c|cccccccl}
      \partial &             -6 & -3 & -2 & -1 & 1 & 2 & 3 & 6  \\
                   \hline
      \# &             5 & 6 & 6 & 381 & 1 & 7 & 2 & 1.
                  \end{array}
\]
\end{Proposition}
\noindent Two aspects of our incomplete numerics are striking.  First, it is somewhat 
surprising that there are at least $408$ number fields with Galois group $\Gamma.2$ and discriminant $2^j 3^k$.  
By way of contrast, the number of fields with Galois group $S_7$, $S_8$,
and $S_9$ and discriminant $\pm 2^j 3^k$ is exactly $10$, at least $72$, and at least $46$ 
respectively \cite{JRGlob}.    Second, the imbalance with respect to $\partial$ 
is quite extreme.    We have not been exhaustive in specializing
our covers and we expect that the $381$ could be increased
somewhat.  By exhaustively specializing Shioda's $W(E_7)^+$ family,
in principle one could obtain the correct
values on the bottom row.  Our expectation however is that
most fields have already been found and so the imbalance
favoring $\Q(i)$ is maintained in the complete numerics.

\subsection{Analysis of ramification} 
\label{analysis} 
 In general, let $K$ be a degree $n$ number field with discriminant $d$ and root 
 discriminant $\delta = |d|^{1/n}$.    It is important to simultaneously consider the Galois 
 closure $K^{\rm gal}$, its discriminant $D$
 and its root discriminant $\Delta = |D|^{1/N}$.  
 For a given field $K$, one has $\delta \leq \Delta$.
 To emphasize the fact that
 the large field $K^{\rm gal}$ is never directly seen in computations, we call $\Delta$ the 
 Galois root discriminant or GRD of $K$.   A GRD $\Delta$ is typically much harder to compute than the corresponding
 root discriminant $\delta$, as it requires good knowledge of higher
 inertia groups at each ramifying prime.   
 
 For sufficiently simple $K$, ramification is thoroughly analyzed by the
 website associated to \cite{JRLoc}, and the GRD $\Delta$ is automatically
 computed.   The 409 fields $K$ contributing to Proposition~\ref{g2fields} 
 are not in the simple range, and we will present one {\em ad hoc}
 computation of a GRD $\Delta$ in the next subsection.  As an illustration
 of the general method, we first consider an easier case here.
 
 For the easier case, take $t=-1$ in \eqref{f28}, which corresponds to 
 $(p,q) = (3,3)$ via $B^*$ and $(a,b) = (-48,9)$ via $B$, both of which
 come from $(u,v) = (1,2)$.   
The discriminant and root discriminant of $K = \Q[x]/f_{28}(-1,x)$ are 
$d =  2^{92} 3^{24}$ and $\delta \approx 25.007$ respectively.  This
root discriminant is much smaller than the minimum $(2^{66} 3^{46})^{1/28} \approx 31.147$ 
appearing in \S\ref{specpi}.   The field $K$ was excluded from consideration 
in \S\ref{specpi} because the Galois group is not $\Gamma.2$
but rather the $336$-element subgroup $PGL_2(7)$.  This
drop in Galois group is confirmed by the factorization of the
resolvent into irreducibles: $f_{36}(-1,x) = x f_{14}(x) f_{21}(x)$.

The group $PGL_2(7)$ can be embedded in $S_8$, which means that
$K^{\rm gal}$ can also be given as the splitting field of a degree eight polynomial.
Such a degree eight polynomial was already found in \cite[Table~8.2]{JRlowgrd}:
\[
f(x) = x^8 - 6 x^4 - 48 x^3 - 72 x^2 - 48 x - 9.
\]
The analysis of ramification is then done automatically by the website associated to
 \cite{JRLoc},
returning for each prime $p$ a slope content symbol $SC_p$ of the form $[s_1,\dots,s_k]_t^u$.  
This means that the decomposition group $D_p$ has order
$p^k t u$, the inertia subgroup $I_p$ has order $p^k t$, and
the wild inertia subgroup $P_p$ has order $p^k$.   
The wild slopes $s_i$ are then rational numbers 
greater than one measuring wildness of ramification, as explained
in \cite[\S3.4]{JRLoc}.

In our $PGL_2(7)$ example, also taking weighted averages to get Galois mean slope \cite[\S3.7]{JRLoc}, 
the result is
\begin{align*}
SC_2 & = [2,3,7/2,9/2]_1^1, & GMS_2 & = \frac{1}{16} \cdot 2 + \frac{1}{8} \cdot 3 + \frac{1}{4}  \cdot \frac{7}{2} + \frac{1}{2}  \cdot \frac{9}{2} = \frac{29}{8}, \\
SC_3 & = [ \;]_7^6, & GMS_3 & = \frac{6}{7}.
\end{align*}
The Galois root discriminant is then $\Delta = 2^{29/8} 3^{6/7} \approx 31.637$.  
This Galois root discriminant is the fifth smallest currently on \cite{JRGlob} 
from a field with Galois group $PGL_2(7)$.  

\subsection{A lightly ramified number field} 
\label{spechard}
Let $K$ be the number field coming from the eight specialization 
points \eqref{repetition}.   Applying {\em Pari}'s \verb@polredabs@ \cite{Pari} to get a canonical
polynomial, this field is $K = \Q[x]/f(x)$ with
 \begin{eqnarray*}
\lefteqn{f(x)=} \\
&& x^{28}-4 x^{27}+18 x^{26}-60 x^{25}+165 x^{24}-420 x^{23}  +798 x^{22}-1440
    x^{21}  +2040 x^{20} \\ &&  -2292 x^{19} +2478 x^{18}-756 x^{17}-657 x^{16}+1464
    x^{15} -4920 x^{14}  +3072 x^{13} \\ && -1068 x^{12}+3768 x^{11}  +1752 x^{10}-4680
    x^9-1116 x^8+672 x^7  +1800 x^6  -240 x^5 \\ && -216 x^4-192 x^3+24 x^2+32 x+4.
\end{eqnarray*}
The field $K$ arises from \eqref{f28} with either $t=4$ or $t = 32/81$, so we also have
its resolvent $K_{36} = \Q[x]/f_{36}(4,x)$ from \eqref{f36}.   Since one of the eight specialization points in 
\eqref{repetition} is $(u,v) = (-4,-3)$, we have also seen this field already 
in the three subsections about $L$-polynomials, \S\ref{UL}, \S\ref{SL}, and \S\ref{GL}.  

Let $K^{\rm gal}$ be the splitting field of $K$.   Calculation of slope contents 
is not automatically done by the website of \cite{JRLoc} because degrees are too large.   The 
proof of the following proposition illustrates 
the types of considerations which are built into
\cite{JRLoc} for smaller degrees.   
\begin{Proposition} 
\label{badanalysis} The decomposition groups of $K^{\rm gal}$ at the ramified primes have invariants as follows:
\begin{align*}
SC_2 & = [2,2,2,3,3]_1^3, & GMS_2 & = \frac{7}{32} \cd 2 + \frac{3}{4} \cd3 = \frac{43}{16}, \\
SC_3 & = [13/8,13/8,11/6]_8^2, & GMS_3 & = \frac{1}{27} \cd \frac{7}{8} + \frac{8}{27} \cd \frac{13}{8}  + \frac{2}{3}  \cd \frac{11}{6}  = \frac{125}{72}. 
\end{align*}
Thus the root discriminant of $K^{\rm gal}$ is $\Delta = 2^{43/16} 3^{125/72} \approx 43.386$.  
\end{Proposition}
\proof The computation is easier at the prime $p=3$ and so we do it first. 
The field $K$ factors $3$-adically as $K_{27} \times \Q_3$ with $K_{27}$ having discriminant
$3^{46}$.  The exponent arises from three slopes $s_1 \leq s_2 \leq s_3$ via
$46 = 2 s_1 + 6 s_2 + 18 s_3$.    One of the two degree $63$ resolvents, computed by {\em Magma},
 factors $3$-adically as $K_{54} \times K_9$, with
 $K_9 \cong  \Q_3[x]/(x^9+6 x^5+6)$ having slope content  $[13/8,13/8]_8^2$.  This
 forces the remaining slope of $K_{27}$ to be $s_3 = 11/6$. 
 The inertia group $D_3$ is thus the maximal subgroup $3_{+}^{1+2}:8:2$ of $\Gamma.2$, with
 slope content $[13/8,13/8,11/6]_8^2$.   

Moving on to the prime $2$, the field $K$ factors $2$-adically as $K_{16} \times K_{12}$.  Here $K_{16}$ is totally ramified
of discriminant $2^{42}$.  The complement $K_{12}$ contains
the unramified cubic extension of $\Q_2$ and has discriminant $2^{24}$.   Since the group $S_{16} \times (S_4 \wr C_3)$ 
does not contain an element of cycle structure either $8^321$ or $8^3 4$, the decomposition
group $D_2$ cannot contain an element of order eight.  So even though 
$\mbox{ord}_2(|\Gamma.2|)=6$, there can be at most five wild slopes.   

The resolvent $K_{36}$ factors as $K_{16} \times K_{12} \times K_8$, with $K_8 \cong \Q_2[x]/(x^8+2x^7+2)$
having discriminant $2^{14}$ and slope content $[2,2,2]_1^3$.   Thus we have found three slopes to be $2$, $2$, and $2$.   If we can find two more wild slopes we will have identified
all wild slopes.   

The field $K_8$ and the sextic field $K_6 = \Q_2[x]/(x^6+x^2+1)$ with discriminant $2^6$
have the same splitting field.   The latter is a subfield of $K_{12}$ showing that 
$(24-6)/6 = 3$ is a fourth $2$-adic slope.   In fact, since both
involutions in $\Gamma.2$ have cycle type $2^{12} 1^4$ and therefore must appear
in the degree $12$ factor, $3$ is the largest wild slope.     

The quartic subfield of $K_8$ is $K_4 = \Q_2[x]/(x^4+2 x^3 + 2 x^2 + 2)$ with discriminant $2^6$. 
Computation shows it is a subfield of $K_{16}$.
 So the remaining slope $s$ satisfies 
$1 \cd 2 + 2 \cd 2 + 4 \cd s + 8 \cd 3 = 42$ and must also be $3$.  
The tame degree $I_2/P_2$ can only be $1$, as the only other possibility
 $t=3$ would force
$u=2$ and $\Gamma.2$ does not contain a solvable subgroup of order a multiple
of $2^6 3^2=576$.     Thus $D_2$ has order $96$
and slope content $[2,2,2,3,3]_1^3$.   \qed

The Galois root discriminant $\Delta \approx 43.386$ is very low,
as is clear from the discussion in \cite{JRlowgrd}, as updated in \cite[Table~9.1]{JRGlob}.   In fact,
the field $K^{\rm gal}$ is  a current record-holder, in the sense that all known 
  Galois fields with smaller root discriminants involve only
simple groups of size smaller than 6048 in their Galois groups.

 \end{document}